\documentclass[12pt]{amsart}

\date{\today}
\usepackage{amsmath,amsthm,amssymb}
\newtheorem{lemma}{Lemma}[section]
\usepackage{mathptmx}
\usepackage{hyperref}
\usepackage{mathrsfs}
\usepackage[T1]{fontenc}
\usepackage{bbm}

\def\angry#1.{\noindent {\bf #1.}}
\def\1{\mathbbm{1}}

\title[Existence and uniqueness of quadrature domains]{Note about the existence and essential uniqueness of quadrature domains}
\author{Hannah Cairns}
\address{Department of Mathematics and Statistics\\\newline
    \hbox{} \hskip 1em Burnside Hall, Room 1005\\\newline
    \hbox{} \hskip 1em 805 Sherbrooke Street West\\\newline
    \hbox{} \hskip 1em Montreal, Quebec, Canada, H3A 0B9}
\email[H.~Cairns]{hannah.abigail.cairns@gmail.com}
\newtheorem{corollary}[lemma]{Corollary}
\newtheorem{theorem}[lemma]{Theorem}

\makeatletter
\newtheoremstyle{smalllemma}
  {6pt}% space before
  {6pt}% space after
  {\narrower \small}% body font
  {}% indent
  {\bfseries}% header font
  {.}% punctuation
  {.5em}% after theorem header
  {}% header specification (empty for default)
\makeatother
\theoremstyle{smalllemma}

\def\setdelta{\mathbin{\text{\footnotesize $\Delta$}}}
\long\def\eat#1{}

\begin{document}
\maketitle
\tableofcontents

\begin{section}{Introduction}
This note is intended to explain the proof of two facts about quadrature domains: first, they are essentially unique if they exist; and second, they do exist for a large class of weight functions.

The proofs roughly follow Sakai~\cite{sakaiobstacle}, especially the proof of Lemma~\ref{sakaislemma} which is precisely the same as Sakai's Lemma 5.1.
Uniqueness depends on well-known properties of Green's function for the Laplacian, for which we refer to M\"orters and Peres~\cite{peres} and Doob~\cite{doob}. The proof of existence uses measure theory, basic properties of harmonic functions, and the fundamental convergence theorem for subharmonic functions from Doob.

Notation. Lebesgue measure is denoted by $\lambda$. If $A$ and $B$ are two sets, then $A \setdelta B := (A \setminus B) \cup (B \setminus A)$. Two sets are {\it essentially equal} if $\lambda(A \setdelta B) = 0$, and $A$ is {\it essentially contained in $B$} if $\lambda(A \setminus B) = 0$.

If $A, B \subseteq \mathbb R^d$, then we set $d(A, B) = \inf\{|x - y|: x \in A, y \in B\}$.

\begin{subsection}{Subharmonic functions}

Let $\Omega \subseteq \mathbb R^d$ be an open set.

Let $B_r(x)$ be the open ball $\{y \in \mathbb R^d: |y - x| < r\}$, and $B_r := B_r(0)$ be the open ball around zero.
If $h: \Omega \to \mathbb R \cup \{\pm\infty\}$ is a locally integrable function, then let $A_h(x; r)$ be the average of the function on the ball $B_r(x)$: $$A_h(x; r) = {1 \over \lambda(B_r)} \int_{B_r(x)} h(y) \,dy.$$
This integral is well-defined for sufficiently small $r$.
We will say that a locally integrable function $h$ is \emph{subharmonic} at a point $x \in \Omega$ if:
\begin{enumerate}
\item[(a)] the function is upper semicontinuous at $x$, so $\limsup_{y \to x} h(y) \le h(x)$,
\item[(b)] and it is bounded above by its averages on small balls around $x$:
$$h(x) \le A_h(x; r) \qquad\text{for all sufficiently small }r.$$
\end{enumerate}
If $h$ is subharmonic at every $x \in E \subseteq \Omega$, then it is \emph{subharmonic on $E$}.
If $-h$ is subharmonic, then $h$ is said to be \emph{superharmonic}.

\begin{subsubsection}{Some properties of subharmonic functions}
\label{subharmonicproperties}
Let $\Omega$ be an open subset of $\mathbb R^d$, and $h$ be subharmonic on $\Omega$. We recall some basic facts.

We have the mean value property: if $B_r(x)$ is any ball with closure contained in $\Omega$, then $h(x)$ is less than or equal to its average on that ball. We also have the maximum principle: if $F \subseteq \Omega$ is compact, then the maximum of $h$ on $F$ is attained at some point on the topological boundary $\partial F$.

If $h$ is twice continuously differentiable on $\Omega$, then $h$ is subharmonic on $\Omega$ if and only if $\nabla^2 h \ge 0$ on $\Omega$.
If $h$ is both subharmonic and superharmonic, then $h$ is said to be \emph{harmonic}, and every harmonic function is smooth, that is, infinitely continuously differentiable.

\medskip
These statements still hold if we replace $\Omega$ by a smaller open set. We will later be able to say something about the Laplacian of a function $h$ that is subharmonic on a general measurable set $E$, in Theorem~\ref{superharmonicmeansnegative}.

\eat{
\end{subsubsection}
\begin{subsubsection}{Limits of radial averages}

We say that a locally integrable function $h$ on $\Omega$ is a \emph{limit of radial averages} on $C$ if, for every point $x \in C$,
$h(x)$ is the limit of $A_h(x; r)$ as $r \to 0$.

A subharmonic function is always a limit of radial averages, because
$$h(x) \le \inf_{0 < r < s} A_h(x; r) \le \limsup_{y \to x} h(y) \le h(x).$$

We call a function \emph{subharmonic on average} on $C$ if it is a limit of radial averages and $h(x) \le A_h(x; r)$ for every sufficiently small $r$.
This condition is strictly weaker than subharmonicity. For example, the sign function is subharmonic on average, but $A_{\text{sign}}(1; 2) = \frac14(-1 + 1 + 1 + 1) = 1/2$.

}

\end{subsubsection}
\end{subsection}

\begin{subsection}{Quadrature domains}
We say a function $w$ is a {\it weight function} if it is bounded, nonnegative, and measurable.
A {\it quadrature domain} for a weight function~$w$ is a bounded open set~$\Omega$ so that we have $w \equiv 0$ outside $\Omega$, and
$$\int hw \,dx \le \int_\Omega h \,dx$$
for every integrable subharmonic function $h$ on $\Omega$.

\medskip
If $h$ is harmonic and bounded on $\Omega$, then $+h$ and $-h$ are subharmonic, so $\int hw \,dx \le \int_\Omega h \,dx$, and $\int (-h)w \,dx \le \int_\Omega (-h) \,dx$, which gives us the opposite inequality $\int_\Omega h \,dx \le \int hw \,dx$. Therefore, $\int hw \,dx = \int_\Omega h \,dx$ in this case.

We can plug in the harmonic functions $h \equiv 1$ and $h = x_i$ to see that the quadrature domain has measure exactly $\int w \,dx$,
and its centre of mass is the same as the centre of mass of the weight function:
$\int_\Omega x_i \,d\lambda = \int x_i w \,d\lambda$, so$${\int_\Omega x_i \,d\lambda \over \lambda(\Omega)} = {\int x_i w \,d\lambda \over \int w \,d\lambda}.$$

If $h = |x|^2 = x_1^2 + \cdots + x_d^2$, then $\nabla^2 h = 2d$, so $h$ is subharmonic, and
$$\int |x|^2 w \,d\lambda \le \int_\Omega |x|^2 \,d\lambda.$$
It follows that the moment of inertia of $\Omega$ is at least as large as the moment of inertia of $w$.
This is consistent with the intuitive description that we get the quadrature domain by spreading out the mass in $w$ in a radially symmetric way.

Quadrature domains are not unique. For example, if we choose $w = 3 \1_{(1, 2)} + 3\1_{(4,5)}$, then $\Omega = (0, 6)$ and $\Omega' = (0, 3) \cup (3, 6)$ are both quadrature domains for $w$.
However, we will see in the next section that quadrature domains for the same weight function are essentially equal.

\def\energy{{\mathscr E}}

\end{subsection}
\end{section}
\begin{section}{Uniqueness of quadrature domains}

We will prove in this section that if a weight function $w$ has two quadrature domains $\Omega$ and $\Omega'$, then $\Omega$ is essentially equal to $\Omega'$.
We need a large supply of subharmonic functions, which come from Green's function.

\begin{subsection}{Green's function}

\begin{subsubsection}{Motivation: What is Green's function?}

Loosely speaking, if we have a differential operator $A = \sum_{\alpha} a_\alpha \partial^\alpha$, we say that a kernel function $g(x, y)$ is ``Green's function'' if the kernel map
$$K_g h(x) = \int h(y) g(x, y) \,dy$$
is a right inverse for $A$. In other words, we have $A K_g h = h$ for any nice $h$.

The function spaces are left deliberately vague.
This is not a precise definition; it is a broad term for a class of similar objects.

We find Green's function for the operator $-\nabla^2$, where
$$\nabla^2 = \sum_i \partial_i \partial_i = \left( {\partial^2 \over \partial x_1^2} + \cdots + {\partial^2 \over \partial x_d^2} \right).$$
We will produce two flavours: an ``unrestricted Green's function,'' and a ``restricted Green's function'' on any bounded open set.

\end{subsubsection}
\begin{subsubsection}{The unrestricted Green's function}
Let $C_d := 2\pi^{d/2} / \Gamma(d/2)$, the area of the unit sphere in $\mathbb R^d$. It turns out that we want our kernel function $G(x,y)$ to be a function $f(|x-y|)$ where $f$ solves $f'(r) = -1/C_d r^{d-1}$.

We make a choice of arbitrary constant and get $$G(x, y) = \begin{cases}\displaystyle{1 \over (d-2) C_d} |x-y|^{2-d}&\text{ if }d \ne 2\\[1em]\displaystyle-{1 \over 2\pi} \log |x-y|&\text{ if }d=2\\[1em]\displaystyle-{1\over2} |x-y|&\text{ if }d=1\end{cases}.$$

If $d = 3$, the constant is chosen so that the function is zero at infinity. Otherwise, we pick a constant that makes the formula simple. In every dimension, this function is smooth on $\mathbb R^d \times \mathbb R^d \setminus \{x = y\}$, and it's symmetric in its arguments, $G(x,y) = G(y,x)$.

Let $G_x(y) := G(x,y)$. If $d = 1$, then the function is locally bounded, and if $d \ge 2$, then it has a singularity but is integrable on any bounded set. One can check that $\nabla^2 G_x = 0$ wherever $G_x$ is smooth, so it's harmonic on the set $\mathbb R^d \setminus \{x\}$. The value at $x$ is $+\infty$, so $G_x$ is superharmonic everywhere.

This gives us a large selection of integrable subharmonic functions to use in the quadrature domain inequality $\int hw \,dx \le \int_\Omega h \,dx$: we can set $h = -G_x$ for any point $x \in \mathbb R^d$, and we can set $h = G_x$ if $x \notin \Omega$.

\end{subsubsection}
\begin{subsubsection}{$G$ is Green's function for $-\nabla^2$}
We check that these kernels we have defined are right inverses of $-\nabla^2$ in the sense of distributions.

Let $\Omega$ be an open subset of $\mathbb R^d$. If $h: \Omega \to \mathbb R$ is smooth and supported on a compact subset of $\Omega$, we say it is a \emph{test function}.

\begin{lemma} \label{test} If $h$ is a test function, then $\int G(x, y) (-\nabla^2 h(y)) \,dy = h(x)$. \end{lemma}

\angry Sketch of proof. We are going to evaluate $\int G_x (-\nabla^2 h) \,dy$. First, if $y \ne x$, then $\nabla^2G_x(y) = 0$, and we can write the integral as a divergence:
\begin{align*}
    G_x (-\nabla^2 h) &= h\nabla^2 G_x- G_x \nabla^2h
    \\&=\nabla \cdot (h \nabla G_x- G_x \nabla h)
\end{align*}

Integrate this divergence over the complement of a small ball $B_\varepsilon(x)$, and then make it an integral over the boundary using Stokes' theorem:
$$\int_{\mathbb R^d \setminus B_\varepsilon(x)} G(x, y) (-\nabla^2 h(y)) \,dy = \int_{|y-x|=\varepsilon} \left(G_x {\partial h \over \partial r} - h {\partial G_x \over \partial r}\right) \,ds.$$

The surface area of the sphere $|y-x|=\varepsilon$ is $C_d \varepsilon^{d-1}$. Green's function is smaller, $G(x, y) = O(1/\varepsilon^{d-2})$ if $d \ne 2$ and $O(\log 1/\varepsilon)$ if $d = 2$; in either case, the integral of $G_x$ over the surface of the sphere is $o(1)$.
On the other hand, the integral of $-h\partial G_x / \partial r = h(y)/C_d \varepsilon^{d-1}$ is going to be $${1 \over C_d \varepsilon^{d-1}} \int_{|y-x| = \varepsilon} h(y)\,dy$$
which is the average of $h$ on the sphere $|y-x| = \varepsilon$. Of course, $h$ is continuous, so that average converges to $h(x)$ as $\varepsilon \to0$. Therefore
\begin{align}
    \int_{\mathbb R^d} G_x(y) (-\nabla^2 h(y)) \,dy
        &= \lim_{\varepsilon \to 0} \int_{\mathbb R^d \setminus B_\varepsilon(x)} G_x(y) (-\nabla^2 h(y)) \,dy\label{greenier}
        \\&=\lim_{\varepsilon \to 0} {1 \over C_d \varepsilon^{d-1}} \int_{|y-x|} h(y) \,dy\notag
        \\&= h(x).\notag
\end{align}

The identity in the first line is true because $G_x$ is a locally integrable function and $h$ is compactly supported, so $G_x \nabla^2 h$ is integrable on $\mathbb R^d$. \qed

\medskip

\begin{corollary}
If $\varphi$ is a distribution, then $-\nabla^2 (G_x * \varphi) = \varphi$.\label{dist}
\end{corollary}

\angry Proof. Let $h$ be a test function. Then by definition, $$-\nabla^2(G_x * \varphi)[h] = (G_x * \varphi)[-\nabla^2 h] = \varphi[G_x * -\nabla^2 h] = \varphi(h).$$
\qed

\end{subsubsection}
\begin{subsubsection}{Green's function of a bounded open set}

\eat{

\angry Proof. Let $h$ be a compactly supported test function. Integrate $G \nabla^2 h = \nabla \cdot (G \nabla h - h \nabla G)$ on $\mathbb R^d \setminus B_\varepsilon(x)$, and then use Stokes's theorem to get \begin{align*}\int_{\mathbb R^d \setminus B_\varepsilon} G(x, y) \nabla^2 h(y) \,dy = \int_{\partial B_\varepsilon(x)} \bigg(G {\partial h \over \partial {\bf n}} - h {\partial G \over \partial {\bf n}}\bigg) \,dS,\end{align*} where $r := ||x - y||$ and ${\bf n}$ is the inner normal of the ball $B_\varepsilon(x)$.

Green's function is a constant times $r^{2-d}$ (or $\log r$ in two dimensions), so the first part of the integrand is
$\int_{\partial B_\varepsilon(x)} G {\partial h \over \partial {\bf n}} \,dS = \int_{\partial B_\varepsilon(x)} O(\varepsilon^{2-d}) = O(\varepsilon)$
for $d \ge 3$ and~$O(\varepsilon \log 1/\varepsilon)$ for $d = 2$. In both cases it vanishes as $\varepsilon \to 0$.

The second part of the integrand is $h \partial G / \partial {\bf n} = -h \partial G / \partial r = 1/(C_d r^{d-1})$, the reciprocal of the surface area of $\partial B_r$. So,
$$\int_{\partial B_\varepsilon(x)} h {\partial G \over \partial {\bf n}} \,dS = {\displaystyle \int_{\partial B_\varepsilon(x)} h(y) \,dS \over \displaystyle \int_{\partial B_\varepsilon(x)} 1 \,dS} = h(x) + O(\varepsilon).$$ Then $\int_{\mathbb R^d \setminus B_\varepsilon} G \nabla^2 h \,dy = O(\varepsilon) - [h(x) + O(\varepsilon)] = -h(x) + O(\varepsilon)$. The integral is finite; take $\varepsilon \to 0$ to get the result $\int_{\mathbb R^d} G \nabla^2 h \,dy = -h(x)$.\qed
}

We define the restricted Green's function in terms of Brownian motion as on page 80, section 3.3 of M\"orters and Peres~\cite{peres}.
We start with a bounded open set $\Omega$ of $\mathbb R^d$.

Let $B_x(t)$ be Brownian motion started at $x \in \mathbb R^d$. It is a random continuous curve in $\mathbb R^d$. Let $T_x := \min\{t \ge 0: B_x(t) \notin \Omega\}$ be the first time that $B_x(t)$ leaves $\Omega$. The \emph{diffusion kernel restricted to $\Omega$} is the continuous function $\mathfrak{p}_\Omega: (0, \infty) \times \Omega \times \Omega \to [0, \infty)$ satisfying
$$\mathbb P\left[B_x(t) \in A \text{ and } t < T_x\right] = \int_A \mathfrak{p}_{\Omega}(t; x, y) \,dy$$ for any measurable set $A$. That is, if we start the walk at $x$ and stop it when it leaves $\Omega$, then at time $t$, the probability density of the walk is $\mathfrak{p}_{\Omega}(t; x, y)$. That function does exist; see section 3.3 of M\"orters and Peres~\cite{peres}.

Set $G_\Omega(x, y) = \int_0^\infty \mathfrak{p}_\Omega(t; x, y) \,dt$. We call this the \emph{restricted Green's function of $\Omega$}. It keeps track of the average amount of time that the walk spends in a region before it leaves the set $\Omega$, in the sense that $$\mathbb E[\text{time that the walk spends in $A$ before leaving $\Omega$}] = \int_A G_\Omega(x, y) \,dy.$$
If $d \ge 3$, Brownian motion eventually wanders off to infinity, and the time spent in any bounded region is finite. In this case, the unrestricted Green's function can be defined in the same way, with $\Omega = \mathbb R^d$: $$\mathbb E[\text{time that the walk spends in $A$}] = \int_A G(x,y)\,dy.$$ That is the connection between the two functions. It's also possible to define $G_\Omega(x,y)$ in terms of a maximization problem as in Section~\ref{maximizationproblem}.

%If we are in three or more dimensions, and if $\Omega$ is the whole space $\mathbb R^d$, then the diffusion kernel is $\mathfrak p(t;x,y) = (2\pi t^2)^{-d/2} e^{-|x-y|^2/2t}$, and it's an interesting exercise to check that $G_\Omega(x,y) = \int_0^\infty \mathfrak p(t;x,y) \,dt$ is $G(x,y)$.

We will use these well-known properties of the restricted Green's function without proof:
\begin{itemize}
\item $G_\Omega(x,y) \ge 0$ and $G_\Omega(x,y) = G_\Omega(y,x)$.
\item If $G_{\Omega,x}$ is the function $y \mapsto G_\Omega(x,y)$, then $G_{\Omega,x}$ is superharmonic on the whole set $\Omega$ and harmonic on $\Omega \setminus \{x\}$. %$y \mapsto G_\Omega(x, y)$ is superharmonic on $\Omega$ and harmonic on $\Omega \setminus \{x\}$.
\item $G_\Omega(x, y) > 0$ when $x,y$ are in the same connected component of $\Omega$.
\item $G_x - G_{\Omega,x}$ is harmonic on $\Omega$.
\end{itemize}

We also want to know that $G_{\Omega,x}$ is integrable on $\Omega$, which we will prove. If $d \ge 3$, then the time that the walk spends in $A$ before $T_x$ is at most the total time we spend in $A$, so $0 \le G_{\Omega,x} \le G_x$ and $\int_\Omega G_{\Omega,x}(y) \,dy \le \int_\Omega G_x(y) \,dy < \infty$.

The situation is more complicated if $d = 2$, but M\"orters and Peres, Lemma 3.37, tells us that $|G(x, y) - G_\Omega(x, y)| \le (1/\pi) \log R/r$, where $r := d(x, \Omega^c)$ and $R := \inf\{r: B_r(x) \supseteq \Omega\}$. This is a finite constant for fixed $x$, so $|G_{\Omega,x}| \le |G_x| + C$ and therefore it is also integrable.
\end{subsubsection}
\end{subsection}
\begin{subsection}{The extended Green's function}
Writers who describe $G_D$ \!as ``the Green's function'' should be condemned to
differentiate the Lebesgue's measure using the Radon-Nikodym's theorem.

\hfill --- Joseph Doob

\medskip

There is a subharmonic extension of $G_{\Omega,x}$ to the set $\mathbb R^d \setminus \{x\}$, and despite the above proscription, we call it {\it the extended Green's function}.

\begin{theorem}If $\Omega$ is a bounded open set with Green's function $G_\Omega$, there is a nonnegative extension  $G_\Omega^e: \Omega \times \mathbb R^d \to [0, \infty]$ so that:
\begin{itemize}
\item $G^e_\Omega(x, y) = G_\Omega(x, y)$ when $y \in \Omega$.
\item $G^e_\Omega(x, -)$ is subharmonic on $\mathbb R^d \setminus \{x\}$ and superharmonic on $\Omega$.
\item $G^e_\Omega(x, -)$ is zero at almost every point in $\mathbb R^d \setminus \Omega$.
\item $G^e_\Omega(x, -)$ is integrable.
\end{itemize}
\end{theorem}

\angry Reference. This is (c) of Doob's Theorem 1.VII.4~\cite{doob}.
Our set $\Omega$ is bounded and $d \ge 2$, so it is Greenian and the theorem applies.

Doob's conclusion is that the function is zero at ``quasi every finite point,'' or in other words on the complement of a polar set. This implies that it is zero on a set of full measure, because polar sets have measure zero.

The extension is integrable because it's equal to $G_\Omega$ on $\Omega$ and zero almost everywhere outside of it.\qed

\medskip
The above theorem is true for \emph{every} bounded open set, even if it has sharp cusps, an infinite number of small holes, or a boundary that has positive Lebesgue measure. This generality is very important in our setting.

\end{subsection}
\begin{subsection}{Monotonicity of quadrature domains}\label{monotonicityofquadraturesets}

\begin{theorem}If $w \le w'$ are two weight functions and $C, D$ are quadrature domains for $w$ and $w'$, then $C$ is essentially contained in $D$.
\label{quadraturesessentiallymonotonic}\end{theorem}

\angry Proof.
Let $E$ be a connected component of $C$.
We will prove that $E \setminus D$ has zero measure.
If it's empty, we are done.
Otherwise, choose $x \in E \setminus D$.
Then $G^e_C$ is subharmonic and integrable on $D \subseteq \mathbb R\setminus\{x\}$, so by the definition of quadrature domains, $\int G^e_C(x, y) w'(y) \,dy \le \int_D G^e_C(x, y) \,dy$.

On the other hand, $-G^e_C(x, -)$ is subharmonic on $C$, and again by the definition of a quadrature domain, $-\int G^e_C(x,y) w\,dy \le -\int_C G^e_C(x,y) \,dy$.

We can chain those inequalities together to get this:
$$\int_C G^e_C(x,y) \,dy \le \int G^e_C(x,y) w(y) \,dy \le \int G^e_C(x, y) w'(y) \,dy \le \int_D G^e_C(x,y) \,dy.$$
Subtract $\int_{C \cap D} G^e_C \,dy$ from both sides to see that $\int_{C \setminus D} G^e_C(x, y) \,dy$ is at most $\int_{D \setminus C} G^e_C(x, y) \,dy$, which is zero because $G^e_C(x, -)$ is zero almost everywhere on the complement of $C$.

Green's function $G_C(x,-)$ is strictly positive on~$E$, so $E \setminus D$ must have zero measure; otherwise $\int_{C \setminus D} G^e_C(x,y) \,dy$ would have been positive. An open set in Euclidean space has only countably many components, so $$C \setminus D = \bigcup_{E\text{ component of }C} E \setminus D$$ also has zero measure, and $C$ is essentially contained in $D$.\qed

\begin{corollary}
Quadrature domains are essentially unique.\label{quadraturesessentiallyunique}
\end{corollary}

\angry Proof. By the lemma, two quadrature domains for the same weight are essentially contained in each other, so their set difference has measure zero.

\end{subsection}
\end{section}
\begin{section}{Positivity of the Laplacian}\label{pos}
\begin{subsection}{Positive distributions}

Let $\Omega$ be an open subset of $\mathbb R^d$, and recall that a {\it distribution} on $\Omega$ is a continuous linear functional on the space of test functions on $\Omega$.

(We recall from distribution theory that a linear functional $\psi$ on the space of test functions is continuous if and only if, for every compact $F \subset \Omega$, the restricted map $\psi|_{C^\infty_c(F)}: C^\infty_c(F) \to \mathbb R$ is continuous in some $\Vert \ \Vert_{C^n(F)}$ norm.
In particular, measures and locally integrable functions are distributions, and the space is closed under differentiation.)

Let $D'(\Omega)$ be the vector space of distributions.
Let $\psi \in D'(\Omega)$.
Then $\partial_i \psi$ is the distribution $h \mapsto -\psi[\partial_i h]$, and $\nabla^2 \psi$ is the distribution $h \mapsto \psi[\nabla^2 h]$.

If $\mu$ is a locally finite measure or a signed measure\footnote{Our signed measures are always bounded, $|\nu|(\Omega) < \infty$.}
we write the corresponding distribution $h \mapsto \int h \,d\mu$ as $d\mu$.
In the same way, if $f$ is locally integrable, we write the distribution $h \mapsto \int hf \,d\lambda$ as $f \,d\lambda$.

(This is a little different from the usual notation, where distributions from measures are written as $\mu$ and distributions from functions are written as $f$.)

A distribution $\psi$ is \emph{positive} if $\psi[h] \ge 0$ whenever $h \ge 0$, and it is \emph{negative} if $\psi[h] \le 0$ whenever $h \ge 0$.
We write $\psi \le \psi'$ if $\psi' - \psi$ is positive.

\begin{subsubsection}{The Laplacian of a subharmonic function}

We will now prove that the distributional Laplacian of a subharmonic function on $\Omega$ is a locally finite measure.

\begin{theorem}\label{superharmonic}
If $f$ is subharmonic on $\Omega$, then $\nabla^2 (f \,d\lambda)$ is positive on $\Omega$.
\end{theorem}

\def\supp{\mathop{\text{supp}}}%
\angry Proof. Fix $h \in C^\infty_c(\Omega)$ with $h \ge 0$. Let $h$ be supported on compact $K \subset \Omega$.
Let $\varphi$ be a positive mollifier, that is, an infinitely differentiable nonnegative function on $\mathbb R^d$ with $\varphi(x) \equiv 0$ for $|x| \ge 1$ and $\int \varphi \,dx = 1$.

Let $\varphi_n(x) := n^d \varphi(nx)$.
Let $n > d(F, \Omega^c)$. Then $f_n := f * \varphi_n$ is defined and infinitely differentiable on a neighbourhood of $K$.
It is also subharmonic:
\begin{align*}f * \varphi_n(x) &= \int_{B_{1/n}} f(x-y) \varphi_n(y) \,dy \\&\le {1 \over \lambda(B_r)} \int_{B_{1/n}} \int_{B_r} f(x-y+z) \varphi_n(y) \,dz\,dy
\\&={1 \over \lambda(B_r)} \int_{B_r} f*\varphi_n(x+z) \,dz.\end{align*}
Therefore, $\nabla^2(f * \varphi_n)$ exists and is nonnegative.

It is well-known that $f * \varphi_n \to f$ in $L^1(F)$. So,
\begin{align*}
\nabla^2 (f\,d\lambda)(h) = \int f \, \nabla^2h \,dx &= \lim_{n \to \infty} \int (f * \varphi_n) \, \nabla^2h \,dx
\\&= \lim_{n \to \infty} \int \nabla^2 (f * \varphi_n) \, h \,dx
\\&\ge 0.
\end{align*}
That is true for every nonnegative $h \in C^\infty_c(\Omega)$, so $\nabla^2 f$ is positive on $\Omega$.\qed

\begin{lemma}If $\psi$ is a positive distribution on $\Omega$, then there is a locally finite measure $\mu$ on $\Omega$ with $\psi(h) = \int h \,d\mu$ for every test function $h \in C^\infty_c(\Omega)$.\label{positivemeansmeasure}\end{lemma}

\angry Proof. See~Rudin~\cite{rudinfunctional}, chapter 6 exercise 4. We sketch the proof.

Let $\psi$ be a positive distribution. If $K \subset\subset \Omega$ is a compact subset, then there is a nonnegative $h_1 \in C^\infty_c(\Omega)$ that is identically $1$ on $K$.
Positivity says $0 \le \psi(h) \le \psi(h_1)$ if $h \in C^\infty_c(K)$ and $0 \le h \le 1$.
Therefore, the restricted distribution $\psi|_{K}: C^\infty_c(K) \to \mathbb R$ is continuous with respect to $\Vert h \Vert_\infty$.

By the Riesz representation theorem and the positivity, there is a finite measure $\mu_K$ with $\psi(h) = \int h \,d\mu_K$ for every test function $h$ supported on the compact set.
These restricted measures are compatible, and we can combine them to get a locally finite measure $\mu$ on $\Omega$ with~$\psi(h) = \int h \,d\mu$ for $h \in C^\infty_c(\Omega)$.

\begin{corollary}If $f$ is subharmonic on $\Omega$, then there exists a locally finite measure $\mu$ on $\Omega$ with $\nabla^2 (f \,d\lambda) = d\mu$.\label{subharmonicmeasure}
\end{corollary}

\angry Proof. The distributional Laplacian $\nabla^2 (f \,d\lambda)$ is positive by Theorem~\ref{superharmonic}, so there is a locally finite measure $\mu$ with $\nabla^2(f \,d\lambda) = d\mu$ by Lemma~\ref{positivemeansmeasure}.\qed

\end{subsubsection}
\begin{subsubsection}{What is coming next}
Corollary~\ref{subharmonicmeasure} tells us a lot about functions which are subharmonic on open sets. What if a function $f$ is subharmonic on a general measurable set?

It turns out that if we already know that the distributional Laplacian $\nabla^2(f \,d\lambda)$ is a signed measure $d\nu$,
 then we can get very precise information: if $f$ is subharmonic on a measurable set $E$, then $\nu$ is positive on $E$.

That is a special case of Theorem~\ref{superharmonicmeansnegative}. In the next few sections we will study the relationship between the spherical average function and the distributional Laplacian, and then finally prove that theorem.

\end{subsubsection}
\end{subsection}
\begin{subsection}{The existence of the spherical average function}

Suppose $f$ is locally integrable, $x$ is a point in $\Omega$, and $0 < r < d(x, \Omega^c)$.
Let the average on the sphere of radius $r$ around $x$ be~$$L_f(x; r) := {1 \over C_d} \int_{|z|=1} f(x + r z) \,dz.$$ Again, $C_d$ is the total surface area of the unit sphere.

If the reader is doubtful about the notation $\int_{|z| = 1} \,dz$, we will see a concrete interpretation in the next section,
Section~\ref{multidimensional}.

\begin{lemma}Let $\Omega$ be an open set, and fix a point $x \in \Omega$.

If $f$ is locally integrable on $\Omega$, then $L_f(x; r)$
is defined for almost every $r < d(x, \Omega^c)$, and $\int_0^s r^{d-1} |L_f(x; r)| \,dr < \infty$ for $s < d(x, \Omega^c)$. 

\label{averagelocallyintegrable}\end{lemma}

\angry Proof.
Let $0 < s < R$. $\overline{B_s(x)}$ is compact, so $\int_{B_s(x)} |f| \,dx < \infty$.

Write this as a double integral in polar coordinates $y = x + rz$, where $r > 0$ and $z$ is a point on the unit sphere:
\begin{align*}\int_{B_s(x)} \big|f(y)\big| \,dy &= \int_0^s \left[\int_{|z| = 1} \big|f(x+r z)\big| \,r ^{d-1}\,dz\right] \,dr.\end{align*}

The left-hand integral is finite, so Tonelli's theorem tells us that the integral in brackets is finite for a.e.~$r \in (0, s)$, and therefore a.e.~$r \in (0, R)$.

Therefore, $L_f(x; r) := \smash{{1 \over C_d} \int_{|z| = 1} f(x+rz) \,dz}$ is well-defined a.e., and
\begin{align*}\int_0^s |L_f(x; r)| \,r^{d-1} \,dr &={1 \over C_d} \int_0^s \left|\int_{|z|=1} f(x+rz)\right| \,r^{d-1} \,dz \,dr\\&\le {1 \over C_d} \int_0^s \int_{|z| = 1} |f(x+rz)| \,r^{d-1} \,dz \,dr\\&<\infty.\end{align*}

\begin{corollary} The function $r \mapsto L_f(x; r)$ is locally integrable on $(0, d(x, \Omega^c))$.
\end{corollary}

\angry Proof. If $0 < t < s < d(x, \Omega^c)$, then $(r/t)^{d-1} \ge 1$ when $r \in [t, s]$, so $$\int_t^s |L_f(x; r)| \,dr \le {1 \over t^{d-1}}\int_t^s |L_f(x; r)| \,r^{d-1}\,dr < \infty.$$
Therefore, $r \mapsto L_f(x; r)$ is integrable on compact subsets of $(0, d(x, \Omega^c))$.
\qed

\begin{corollary} \label{a}The average of a locally integrable function $f$ on $B_s(x)$ is
$$A_f(x; s) = \int_0^s {dr^{d-1} \over s^d} L_f(x; r) \,dr.$$

\end{corollary}
\angry Proof. We use the change of variables $y = x + rz$ again to write
\begin{align*}
    \int_{B_s(x)} f(y) \,dy &= \int_0^s \int_{|z| = 1} f(x + rz) \,r^{d-1} \,dr \,dz
    \\&=\int_0^s C_d L_f(x; r) \,r^{d-1} \,dr.
\end{align*}
Dividing this by the integral $\int_{B_s(x)} 1 \,dy = C_d s^d / d$ gives the result.
\qed

\medskip
At this point we know that $L_f(x; r)$ makes sense and is integrable.
The next step is to show that spherical averages are related to the distributional Laplacian by an integral equality. We do that in the next section.

\end{subsection}
\begin{subsection}{A digression: multidimensional polar coordinates}
\label{multidimensional}In case one finds the ``polar coordinates'' $y = x + rz$ above to be a little suspicious, we will provide a concrete interpretation.

We define multidimensional polar coordinates $r, \varphi_1, \ldots, \varphi_{d-2}, \theta$, where
\eat{
\begin{align*}
y_1 &= x_1 + r\cos \varphi_1 \cdots \cos \varphi_{d-3} \cos \varphi_{d-2} \cos \theta,\\
y_2 &= x_2 + r\cos \varphi_1 \cdots \cos \varphi_{d-3} \cos \varphi_{d-2} \sin \theta,\\
y_3 &= x_3 + r\cos \varphi_1 \cdots \cos \varphi_{d-3} \sin \varphi_{d-2},\\
&\ \,\vdots\\
y_{d-2} &= x_{d-2} + r\cos \varphi_1 \cos \varphi_2 \sin \varphi_3,\\
y_{d-1} &= x_{d-1} + r\cos \varphi_1 \sin \varphi_2,\\
y_{d} &= x_d + r\sin \varphi_1.
\end{align*}
}
\begin{align*}
y_1 = x_1 &+ r \sin \varphi_1,\\
y_2 = x_2 &+ r \cos \varphi_1 \sin \varphi_2,\\
y_3 = x_3 &+ r \cos \varphi_1 \cos \varphi_2 \sin \varphi_3,\\
\vdots\\
y_{d-2} = x_{d-2} &+ r \cos \varphi_1 \cdots \cos \varphi_{d-3} \sin \varphi_{d-2},\\
y_{d-1} = x_{d-1} &+ r \cos \varphi_1 \cdots \cos \varphi_{d-3} \cos \varphi_{d-2} \sin \theta,\\
y_{d} = x_{d} &+ r \cos \varphi_1 \cdots \cos \varphi_{d-3} \cos \varphi_{d-2} \cos \theta
\end{align*}

The bounds are $r > 0$, $-{\pi \over 2} \le \varphi_1, \ldots, \varphi_{d-2} \le {\pi \over 2}$, $0 \le \theta < 2\pi$.
They are the usual polar coordinates when $d = 2, 3$. For example, in $d = 3$,
\begin{align*}
y_1 &= x_1 + r \cos \varphi_1 \cos \theta,\\
y_2 &= x_2 + r \cos \varphi_1 \sin \theta,\\
y_3 &= x_3 + r \sin \varphi_1.
\end{align*}

We want to find the determinant of the Jacobian matrix $J$. One way to do this is to calculate the metric $g = JJ^T$, which is diagonal with $g_{rr} = {\partial y \over \partial r} \cdot {\partial y \over \partial r} = 1$, $g_{\varphi_j \varphi_j} = {\partial y \over \partial \varphi_j} \cdot {\partial y \over \partial \varphi_j} = r^2 \cos^2 \varphi_1 \cdots \cos^2 \varphi_{j-1}$, and $g_{\theta \theta} = {\partial y \over \partial \theta} \cdot {\partial y \over \partial \theta} = r^2 \cos^2 \varphi \cdots \cos^2 \varphi_{d-2}$.
We have $\det g = (\det J)^2$, so we take the square root:
$$|\det J| = \sqrt{\det g} = \sqrt{\prod g_{ii}} = r \times (r \cos \varphi_1) \times \cdots \times (r \cos \varphi_1 \cdots \cos \varphi_{d-2}),$$
or $|\det J| = r^{d-1} \cos^{d-2} \varphi_1 \cdots \cos \varphi_{d-2}$.

We don't need to know the sign, but if we want it, we can get it by looking at the point $\varphi_1 = \cdots = \varphi_{d-2} = 0,$ $\theta = 0$, $r = 1$, where $J$ is the identity matrix.
The sign of the determinant is positive there, and $\{\det J \ne 0\}$ is connected and dense, so the sign is nonnegative everywhere.

\medskip
The change-of-variables formula tells us that
\begin{align*}
\int_{\mathbb R^d} f(y) \,dy &= \int_0^\infty \int_{-\pi/2}^{\pi/2} \cdots \int_{-\pi/2}^{\pi/2} \int_0^{2\pi} f(x + rz) \,r^{d-1} \\&\qquad\qquad\qquad\qquad\cos^{d-2} \varphi_1 \cdots \cos \varphi_{d-2} \,d\varphi_1 \cdots \,d\varphi_{d-2} \,d\theta \,dr.
\end{align*}
Now we ask the reader to interpret $\int_{|z| = 1}$ as shorthand for integration over all the coordinates except $r$,
and $dz$ as an abbreviation for the expression $\cos^{d-2} \varphi_1 \cdots \cos \varphi_{d-2} \,d\varphi_1 \cdots \,d\varphi_{d-2} \,d\theta$.

Then we do have the identity $\int f(y) \,dy = \int_0^\infty \int_{|z|=1} f(x + rz) \, r^{d-1} \,dz \,dr$, and the derivation in the last section makes sense.

\begin{subsubsection}{Exercise: the value of $C_d$}
This gives us another expression for $C_d$:
\begin{align*}
C_d
    = \int_{|z| = 1} \,dz
    &= \int_{-\pi/2}^{\pi/2} \cdots \int_{-\pi/2}^{\pi/2} \int_0^{2\pi} \cos^{d-2} \varphi_1 \cdots \cos \varphi_{d-2} \,d\varphi_1 \cdots d\varphi_{d-2} \,d\theta
    \\&=2\pi\prod_{j=1}^{d-2} \int_{-\pi/2}^{\pi/2} \cos^j \varphi \,d\varphi.
\end{align*}
Use the beta integral $\int_0^1 x^{\alpha - 1} (1-x)^{\beta - 1} \,dx = \Gamma(\alpha) \Gamma(\beta) / \Gamma(\alpha + \beta)$ and the special value $\Gamma(1/2) = \sqrt{\pi}$ to show that this is equal to $2\pi^{d/2} / \Gamma(d/2)$.
Hint: evaluate $\int_{-\pi/2}^{\pi/2} \cos^j \varphi \,d\varphi = \Gamma(\frac12) \Gamma(j/2 + \frac12) / \Gamma(j/2 + 1)$.

\end{subsubsection}
\end{subsection}
\begin{subsection}{Spherical averages and the distributional Laplacian}
We can get many weighted integrals of the spherical averages $L_f(x; r)$ by evaluating the Laplacian $\nabla^2(f \,d\lambda)$ on certain nonnegative functions.

\begin{lemma} Let $x \in \Omega$. Let $R := d(x, \Omega^c)$. Let $\eta \in C^\infty_c(0, R)$ with $\eta \ge 0$. Then there is a nonnegative function $h \in C^\infty_c(\Omega)$ with
$$\smash{\nabla^2 (f \,d\lambda)(h) = -\int_0^R \eta'(r) L_f(x; r) \,dr}$$ for every locally integrable function $f$ on $\Omega$.

\label{etalemma}
\end{lemma}

\angry Proof. Let $h(y) := H(|y-x|)$ for $y \in B_R(x)$, where
$$H(r) := \int_{r}^R {\eta(\rho) \over C_d \rho^{d-1}} \,d\rho.$$
The compactly supported function $\eta$ is zero on a neighbourhood of $0$ and $R$,
so $H(r)$ is smooth, constant near $0$, and zero on a neighbourhood of $R$.
Therefore, $h$ is smooth even at $x$, and compactly supported in $B_R(x)$.

If a smooth function $h$ is radially symmetric around a point $x$, then there is a formula for its Laplacian which holds in any $\mathbb R^d$:
$$\nabla^2 h = {1 \over r^{d-1}} {\partial \over \partial r} \bigg[r^{d-1} {\partial \over \partial r} h\bigg]$$where we are using the polar coordinates $y = x + rz$.%
\footnote{
We can get this from the Voss-Weyl formula $\nabla^2 h = {1 \over \sqrt{\det g}} \sum_{ij} {\partial \over \partial \xi_i} (\sqrt{\det g} \,
g^{ij} {\partial \over \partial \xi_j} h)$.
Here $g^{ij} = (g^{-1})_{ij}$. Using the coordinates from Section \ref{multidimensional}, we get the result.
}

The radial derivative of $h$ is of course ${\partial h \over \partial r} = H'(r) = -\eta(r) / C_dr^{d-1}$, which means the Laplacian is $\nabla^2 h = -\eta'(r) / C_dr^{d-1}$.
Then
\begin{align*}
\nabla^2(f \,d\lambda)(h) &= \int_{B_R(x)} f \,\nabla^2 h \,dy
\\&=-\int_{B_R(x)} f(y) {\eta'(r) \over C_d r^{d-1}} \,dy \\&=-\int_0^R \int_{|z|=1} f(x+rz)\,  {\eta'(r) \over C_d r^{d-1}} \,r^{d-1} \,dz \,dr
\\&=-\int_0^R {\eta'(r) \over C_d} \int_{|z| = 1} f(x + rz) \,dz \,dr
\\&=-\int_0^R \eta'(r) L_f(x; r) \,dr.\end{align*}
This is the result.\qed

\medskip

We will use this basic result to evaluate the differences $A_f(x; s) - A_f(x; t)$
in terms of the distributional Laplacian. We will get an especially precise result when $\nabla^2(f \,d\lambda) = d\nu$ for some signed measure $\nu$.

\end{subsection}
\begin{subsection}{Difference of averages: choosing functions for approximation}

\eat{
If $\eta \in C^\infty_c(0, R)$, $\eta \ge 0$, then $\eta'(r)$ is a signed weight function that's smooth and compactly supported on $(0, R)$, and has total weight zero.

The positive weight is left of the negative weight, in the sense that $\int_0^r \eta'(\rho) \,d\rho \ge 0$ and $\int_r^\infty \eta'(\rho) \,d\rho \le 0$.

The converse is also true: any function like that will give us a suitable $\eta \in C^\infty_c(0, R)$, $\eta \ge 0$ by integration.
}

Recall that $A_f(x; r)$ is the average of $f$ on the ball $B_x(r)$.

Let $x \in \Omega$, $0 < r < s < d(x, \Omega^c)$. We construct some nonnegative functions $\eta_m$ that are suitable for Lemma~\ref{etalemma}, and then use it to prove that \begin{equation}\label{theconcentricballsidentity}A_f(x; s) - A_f(x; t) = \lim_{m \to \infty} \nabla^2(f \,d\lambda)(h_m)\end{equation}
where $h_m(y) := \int_{|y-x|}^\infty \eta_m(r) / C_d r^{d-1} \,dr$ as in the lemma.

Write the formula in Corollary \ref{a} as $A_f(x; s) = \int_0^\infty \1_{r < s} {dr^{d-1} \over s^d} L_f(x; r) \,dr$.
Then we can write the difference $A_f(x; t) - A_f(x; s)$ as
$$A_f(x; t) - A_f(x; s) = \int_0^\infty \left[\1_{r < t} {dr^{d-1} \over t^d} - \1_{r < s} {dr^{d-1} \over s^d}\right] L_f(x; r) \,dr.$$
We want to connect this to Lemma~\ref{etalemma}. Let $w$ be the function \begin{equation}\label{definitionofw}w(r) := \begin{cases}dr^{d-1}&\text{if }r<1,\text{ and}\\0&\text{if }r\ge1.\end{cases}\end{equation}
Let $W(r) := \min\{1, r^d\}$ be the integral of $w$ from $0$ to $r$.
Then the expression in brackets is the derivative of $W(r/t) - W(r/s)$, where it is differentiable.
Unfortunately, $W$ isn't smooth, so we can't plug it into Lemma \ref{etalemma} directly and we have to approximate.

Let $w_m$ be a sequence of compactly supported, nonnegative smooth functions $w_m \in C^{\infty}_c(0, \infty)$ which increase to $w$. Let $W_m(r) := \int_0^r w_m(\rho) \,d\rho$.

\begin{lemma}
Given $0<t<s$, the sequence of functions $$\eta_{m, s, t}(r) := W_m(r/t) - W_m(r/s)$$% there are functions $\eta_m \in C^\infty_c(0, \infty)$
converges uniformly to $W(r/t) - W(r/s)$, and each function satisfies $$0 \le \eta_{m,s,t}(r) \le \begin{cases}\min\{1, r^d/t^d\}&\text{if }r < s \\0&\text{otherwise}\end{cases}$$ and $$|\eta_{m,s,t}'(r)| \le \begin{cases}dr^{d-1} / t^d&\text{if }r < s\\0&\text{otherwise}.\end{cases}$$
\label{etamfunctions}\end{lemma}

\angry Proof. Here $\eta_{m,s,t}(r) \ge 0$ because $W_m$ is increasing and $r/t > r/s$, and $\eta_{m,s,t}(r) \le W_m(r/t) \le W(r/t) = \min\{1, (r/t)^d\}$.
The rest is easy.\qed

\eat{Let $w(r) = dr^{d-1}$ for $r < 1$ and $0$ for $r \ge 1$. Let $W(r)$ be the indefinite integral $\int_0^r w(\rho) \,d\rho = \min\{1, r^d\}$.
Then $\chi(r) = W(r/t) - W(r/s)$.

Let $w_m \in C^\infty_c(0, \infty)$ be nonnegative smooth functions with supports contained in $(0, 1)$ that increase pointwise to~$w$. Set $W_m(r) := \int_0^r w_m(\rho) \,d\rho$.

Finally, let~$\eta_m(r) := W_m(r/t) - W_m(r/s)$. These are the approximations. They are nonnegative because $W_m$ is increasing, and they're zero for $r \ge s$ because $W_m(r)$ is constant for $r \ge 1$.
The first bound holds: $$0 \le \eta_m(r) \le W_m(r/t) \le W(r/t) = \min\bigg\{1, \frac{r^d}{t^d}\bigg\}.$$
The bound on the derivative holds:
$$|\eta_m'(r)| = \left|{w_m(r/t) \over t} - {w_m(r/s) \over s}\right| \le \max\left\{{w(r/t) \over t}, {w(r/s) \over s}\right\} \le {dr^{d-1} \over t^d}.$$% and $\eta_m(r) \to \eta(r)$
And it's easy to see that $|W(r) - W_m(r)| \le \int_0^R (w(r) - w_m(r)) \,dr$ converges to zero uniformly in $r$, so $\eta_m \to \chi$ uniformly.\qed}

\begin{lemma}Suppose \label{h}$x \in \Omega$ and $0 < t < s$ with $B_s(x) \subseteq \Omega$. If we choose functions $\eta_{m,s,t}$ as above, then we will have \begin{equation*}\tag{\ref{theconcentricballsidentity}}A_f(x; s) - A_f(x; t) = \lim_{m \to \infty} \nabla^2(f \,d\lambda)(h_{m,s,t,x})\end{equation*}for $h_{m,s,t,x}(y) := \int_{|y-x|}^\infty \eta_{m,s,t}(r) / C_d r^{d-1} \,dr$ as in Lemma~\ref{etalemma}.\end{lemma}

\angry Proof. We start from Corollary \ref{a}, and remember that the functions $w_m(r)$ are functions that increase to $dr^{d-1}$ for $r < 1$.
\begin{align*}A_f(x; s) &= \int_0^s {d r^{d-1} \over s^d}\, L_f(x; r) \,dr \notag \\&= \int_0^\infty \lim_{m \to \infty} {w_m(r/s) \over s} \, L_f(x; r) \,dr.
\end{align*}
That is bounded by $\1_{r < s} (d/s) |L_f(x; r)|$, which is integrable by Lemma \ref{averagelocallyintegrable}.
So we can use dominated convergence to get:
$$A_f(x; s) = \lim_{m \to \infty} \int_0^\infty {w_m(r/s) \over s}\, L_f(x; r) \,dr.$$
Replace $s$ by $t$ and subtract:
\begin{align*}
A_f(x; s) - A_f(x; t) &= \lim_{m \to \infty} \int_0^\infty \left[{w_m(r/s) \over s} - {w_m(r/t) \over t}\right] L_f(x; r) \,dr\notag
\\&=\lim_{m \to \infty} \int_0^\infty (-\eta_{m,s,t}'(r)) L_f(x; r) \,dr\\&=\lim_{m \to \infty} \nabla^2(f \,d\lambda)(h_{m,s,t}).%\label{ballaverage}
\end{align*}
The last equality comes from Lemma~\ref{etalemma}.
This is the desired identity.
\qed

\medskip

If $\nabla^2(f \,d\lambda)$ is a signed measure $d\nu$, we can use this lemma to write the difference~$A_f(x; t) - A_f(x; s)$ as an integral with respect to the signed measure.
\end{subsection}
\begin{subsection}{Difference of averages: a formula for signed measures}

\begin{theorem}
\label{ballcomparison} Suppose $f$ is locally integrable on $\Omega$ and $\nabla^2(f \,d\lambda) = d\nu$ where~$\nu$ is a signed measure.
Let $x \in \Omega$ and $0 < t < s < R = d(x, \Omega^c)$.

Let $$h_{s,t,x}(y) := \int_{|y-x|}^\infty {W(r/t) - W(r/s) \over C_dr^{d-1}} \,dr.$$
Then $\int h_{s,t,x} \,d\nu = A_f(x; s) - A_f(x; t)$.
\end{theorem}

\angry Proof. Let $h_{m,s,t,x}$ be the functions provided by Lemma~\ref{h} with $$A_f(x; s) - A_f(x; t) = \lim_{m \to \infty} \nabla^2(f \,d\lambda)(h_{m,s,t,x}).$$

If we can prove that $\lim_m \int h_{m,s,t,x} \,d\nu = \int h_{s,t,x} \,d\nu$, we will be done.
We can write the difference as $$h_{s,t,x}(y) - h_{m,s,t,x}(y) = \int_{|y-x|}^\infty {W(r/t) - W(r/s) - \eta_{m,s,t}(r) \over C_d r^{d-1}} \,dr.$$
This difference is uniformly bounded in absolute value by
$$\int_{0}^\infty {|W(r/t) - W(r/s) - \eta_{m,s,t}(r)| \over C_d r^{d-1}} \,dr.$$

Lemma~\ref{etamfunctions} says that the integrand is uniformly bounded by $\1_{r < s} r C_d/t^d$, and that it goes to zero pointwise. By the dominated convergence theorem, $\max |h_{s,t,x} - h_{m,s,t,x}| \to 0$, so we do have $\lim_m \int h_{m,s,t,x} \,d\nu = \int h_{s,t,x} \,d\nu$.
\qed\medskip

This lemma will allow us to get an estimate on the difference of averages from weak estimates on the Laplacian.
To get it, we need to know $\int h \,dx$.

\begin{lemma}\label{theintegralofh}Let $x \in \Omega$ and $0 < t < s$. With $h_{s,t,x}$ defined as in Theorem~\ref{ballcomparison}, $$\displaystyle\int h_{s, t,x}(y) \,dy = {1 \over 2(d+2)} (s^2 - t^2).$$\end{lemma}

\angry Proof. We know what the function is, so the proof is a calculation. First,
\begin{align*}\int_{\mathbb R^d} h_{s,t,x}(y) \,dy &= \int_{\mathbb R^d} \left[\int_{|x-y|}^s {W(\rho/t) - W(\rho/s) \over C_d \rho^{d-1}}\,d\rho\right] \,dy
\\&= \int_0^\infty \left[\int_{r}^s {W(\rho/t) - W(\rho/s) \over C_d \rho^{d-1}} \, \,d\rho\right] C_dr^{d-1} \,dr
\\&= \int_0^s {W(\rho/t) - W(\rho/s) \over \rho^{d-1}} \left[\int_0^\rho r^{d-1} \,dr\right] \,d\rho
\\&= \int_0^s (W(\rho/t) - W(\rho/s)) {\rho \over d} \,d\rho.
\intertext{We have $W(r) = \min\{1, r^d\}$, so}
\int_0^s (W(\rho/t) - W(\rho/s)) {\rho \over d} \,d\rho&=\left[\int_0^t {\rho^{d+1} \over d t^d} \,d\rho + \int_t^s {\rho \over d} \,d\rho \right] - \int_0^s {\rho^{d+1} \over d s^d} \,d\rho\\&={t^{2} \over d(d+2) } + {s^2 - t^2 \over 2d} - {s^2 \over d(d+2) }
\\&={1 \over 2(d+2)} (s^2 - t^2).\end{align*}
\vskip-1em\qed

\medskip
We will use this in Lemma~\ref{quadraticatazero} and Theorem~\ref{zeroonboundary} to get quadratic bounds on $A_f(x; s) - A_f(x; t)$ in the case where $\nabla^2(f \,d\lambda) = \rho \,d\lambda$ with $0 \le \rho \le 1$.

\eat{

Suppose $\nabla^2 (f \,d\lambda)$ is a signed measure $d\nu$ where $|\nu| \le C\lambda$. Then
$$|A_f(x; s) - A_f(x; t)| \le \left|\int h \,d\nu\right| \le C \int h \,d\lambda = {C \over 2(d+2)} (s^2 - t^2),$$
so $\lim_{t \to 0} A_f(x; t)$ exists for every $x\in\Omega$. Let the limit be $\bar f(x)$. Take $t \to 0$ above to get the inequality $|A_f(x; s) - \bar f(x)| \le Cs^2/2(d+2)$, which is not only a uniform bound, it's \!{\it quadratic} in the radius \!$s$, just as strong as if we were in one dimension and knew that $f \in C^2$ and $|f''| \le C$.%
\footnote{Exercise: prove that we do really have $|{1 \over 2s} \int_{-s}^s f(y) \,dy - f(0)| \le Cs^2/6$ if $|f''| \le C$.}

}

\end{subsection}
\begin{subsection}{Limits of radial averages}\label{lora}

\medskip

A function is a {\it limit of radial averages at $x$} if it is integrable in a neighbourhood of $x$ and the limit $\lim_{r \to 0} A_f(x; r)$ exists and is equal to $f(x)$. This is strictly weaker than continuity at a point.

A subharmonic function is always a limit of radial averages, because
$$h(x) \le \inf_{0 < r < s} A_h(x; r) \le \limsup_{y \to x} h(y) \le h(x).$$

\begin{subsubsection}{Subharmonicity on average}

We say that a function $f$ is {\it subharmonic on average} at $x$ if it is a limit of radial averages at $x$ and satisfies condition (b) in the definition of subharmonicity. That is, there exists some small $\varepsilon > 0$ so that~$$\lim_{r \to 0} A_f(x; r) = f(x) = \inf_{0 < r < \varepsilon} A_f(x; r).$$

\def\sign{\mathop{\text{sign}}}%
This is strictly weaker than subharmonicity. For example, the sign function is subharmonic on average everywhere,
but its average on the interval $B_2(1) = (-1, 3)$ is $A_f(1; 2) = \frac12$, which is strictly less than $\sign 1 = 1$.

From our point of view, the problem is that the distributional Laplacian of the sign function is too irregular:
it is not a signed measure.\footnote{If $h$ is a test function, then $\smash{\nabla^2(\sign)(h) = \int_{-\infty}^\infty \sign x \,h'' \,dx = -2h'(0)}$. Suppose there were a signed measure $\nu$ with $\int h \,d\nu = -2h'(0)$. Then there would be a constant $C = |\nu|(\Omega) / 2$ with $|h'(0)| \le C \max |h|$ for every test function $h$, but this is absurd.}

\end{subsubsection}
\begin{subsubsection}{...implies positivity of the Laplacian}

We will show that, if the Laplacian is a signed measure $d\nu$, then~$\nu$ is positive on any measurable set where $f$ is subharmonic on average.

We need a consequence of the Lebesgue-Besicovitch theorem.

\begin{lemma}\label{frogfractions} If $E \subseteq \Omega$ is measurable and $\nu$ is a signed measure with~$\nu(E) < 0$, then there is a point $x \in E$ with $$\limsup_{t \to 0} {\nu(B_t(x)) \over \lambda(B_t(x))} < 0.$$\end{lemma}

\angry Proof. Let $\mu = |\nu| + \lambda$, and $f = d\nu / d\mu$. Then $\int f \,d\mu = \nu(E) < 0$, so the set of points where $f$ is negative must have positive measure.% we must have $\mu(\{x: f(x) < 0\}) > 0$.% on a set of $\mu$-positive measure.

By the Lebesgue-Besicovitch differentiation theorem,$$\lim_{t \to 0} {\nu(B_t(x)) \over \mu(B_t(x))} = f(x)$$ except on a set $N$ with $\mu(N) = 0$.
The set of points with $f(x) < 0$ has positive measure, so there must be some point $x$ with $f(x) < 0$ and $x \notin N$.

By definition of the limit, $\nu(B_t(x))$ is negative for small $t$, and $0 \le \lambda \le \mu$, so $1/\lambda(B_t(x)) \ge 1/\mu(B_t(x))$ and~$\nu(B_t(x)) / \lambda(B_t(x)) \le \nu(B_t(x)) / \mu(B_t(x))$.

We therefore have the strict inequality
$$\limsup_{t \to 0} {\nu(B_t(x)) \over \lambda(B_t(x))} \le \limsup_{t \to 0} {\nu(B_t(x)) \over \mu(B_t(x))} = f(x) < 0.$$
That proves the result.
\qed

\begin{theorem} Suppose $f$ is locally integrable on $\Omega$, and $\nabla^2 (f \,d\lambda) = d\nu$ where~$\nu$ is a signed measure.\label{superharmonicmeansnegative}

If $f$ is subharmonic on average on a measurable set $E$, then $E$ is a positive set for~$\nu$, i.e.~$\nu(E') \ge 0$ for $E' \subseteq E$.
\end{theorem}

\angry Proof. Suppose $E$ is not positive. Let $E'$ be a measurable subset of $E$ with negative $\nu$-measure. By Lemma~\ref{frogfractions}, $\exists x \in E'$ with $$\limsup_{t \to 0} {\nu(B_t(x)) \over \lambda(B_t(x))} =-c < 0.$$

Let $s > 0$ be small enough that~$\nu(B_t(x)) / \lambda(B_t(x)) < -c/2$ for $t < s$ and the subharmonic inequality holds for $B_s(x)$.

Theorem~\ref{ballcomparison} and Corollary~\ref{theintegralofh} tell us that, for each $x \in \Omega$ and $0 < t < s$, there is a nonnegative radially symmetric continuous function $h_{s,t,x}$ with
\begin{equation*}\int h_{s, t,x}(y) \,d\nu(y) = A_f(x; s) - A_f(x; t).\end{equation*}

The reader can check from the definition that $h_{s, t, x}(y)$ is supported on $B_s$, radially symmetric, and decreases as $y$ gets farther from $x$.
So if $\alpha > 0$, then $\{h_{s, t, x} > \alpha\}$ is a ball around $x$ of radius less than $s$, and
\begin{align*}\int h_{s, t,x}\,d\nu&=\int_0^\infty\nu(\{h_{s,t,x} > \alpha\}) \,d\alpha
\\&\le -{c \over 2}\int_0^\infty \lambda(\{h_{s, t,x} > \alpha\}) \,d\alpha \\&= -{c \over 2} \int h_{s, t, x} \,d\lambda = -{c \over 4(d+2)} (s^2 - t^2).\end{align*}
The last step is Lemma~\ref{theintegralofh}.

Fix $s$, take $t \to 0$, and use the fact that~$A_f(x; t) \to f(x)$ as $t \to 0$, because~$f$ is a limit of radial averages at $x$.
We get the impossible inequality:
$$0 \le A_f(x;s)-f(x)=\limsup_{t\to0} \int h_{s,t,x} \,d\nu \le -{c \over 4(d+2)}s^2 < 0.$$
So, there is no measurable subset $E'$ with $\nu(E') < 0$.
\qed

\end{subsubsection}
\end{subsection}
\begin{subsection}{On an open set, positivity implies subharmonicity}

Now we will go from the Laplacian to full subharmonicity.

\begin{lemma} Let $f$ be locally integrable on an open set $\Omega$. Suppose $\nabla^2 (f \,d\lambda)$ is positive on $\Omega$. Then there is a subharmonic $\bar{f}$ on $\Omega$ with $\bar f = f$ \!a.e.\!~on $\Omega$. \label{negativemeanssuperharmonic}\end{lemma}

\angry Proof.
Let $x \in \Omega$ and $t < s < d(x, \Omega^c)$.

Theorem~\ref{ballcomparison} tells us that $A_f(x; s) - A_f(x; t) = \nabla^2(f \,dx)(h_{s, t}) \ge 0$ for a certain function $h_{s, t}$, so $A_f(x; t)$ decreases to a limit (possibly $-\infty$) as $t \to 0$.

Let $\bar{f}(x)$ be that limit: $$\bar{f}(x) := \lim_{t \to 0} A_f(x; t) = \inf_{t > 0} A_f(x; t).$$
The Lebesgue differentiation theorem tells us that $A_f(x; t) \to f(x)$ for almost every $x$, so $f = \bar f$ for almost every $x \in \Omega$.

We claim $\bar{f}$ is subharmonic on $\Omega$. It is less than or equal to its averages on balls because $\bar f(x) \le A_f(x; t) = A_{\bar f}(x; t)$,
 so (b) in the definition of subharmonicity is satisfied.
We must prove that $\bar f$ is upper semicontinuous.

Let $x_n$ be a sequence of points in $\Omega$ that converge to $x$.
Let $0 < r < d(x, \Omega^c)$. We can break down $\bar{f}(x_n)$ in the following way:
\begin{align*}\bar{f}(x_n) &= \Big[\bar{f}(x_n) - A_f(x_n; r)\Big]+ \Big[A_f(x_n; r) - A_f(x; r)\Big] +A_f(x; r).\end{align*}
The first summand is nonpositive by definition. The second one converges to~$0$ as~$n \to \infty$, because $A_f(x; r)$ is continuous in $x$. This is a general property of convolutions, but it can be proven directly in this case by writing $$A_f(r; x_n) - A_f(r; x) = {1 \over \lambda(B_r)} \int (\1_{B_r(x_n)} - \1_{B_r(x)}) f(y) \,dy.$$ The integrand is dominated by $|f|$ and converges pointwise to zero.
Take the lim sup of both sides of the equality as $n \to \infty$:
$$\limsup \bar f(x_n) = (\text{nonpositive}) + A_f(x; r)$$ and then take $r \to 0$ to get $\limsup \bar{f}(x_n) \le \bar{f}(x)$. So $\bar{f}$ is upper semicontinuous, and therefore subharmonic.
\qed

\medskip

\angry Note. Theorems~\ref{superharmonicmeansnegative} and~\ref{negativemeanssuperharmonic} together tell us that if $f$ is subharmonic on average, and its Laplacian is a signed measure, then it is subharmonic.

\end{subsection}
\begin{subsection}{The measure of the Laplacian on the zero set}

In what follows, we suppose $f$ is a limit of radial averages, $f \ge 0$, and the distributional Laplacian $\nabla^2(f \,d\lambda)$ is $\rho \,d\lambda$ with $|\rho| \le C$.

We will show that~$\rho = 0$ a.e.~on the zero set $\{x: f(x) = 0\}$, which we will need later for the crucial Corollary~\ref{findthelaplacian}.
First, we show that $f(y)$ converges uniformly to zero as $y$ approaches the zero set.

\begin{lemma} Suppose $f \ge 0$ is a limit of radial averages and $\nabla^2(f \,d\lambda) = \rho \,d\lambda$ with $|\rho| \le C$.
If $f(x) = 0$, then~$f(y) \le 2^dC|y-x|^2$ if $|y-x| < \frac12 d(x, \Omega^c)$.\label{quadraticatazero}
\end{lemma}

By Theorem \ref{ballcomparison}, if $x$ is a point in $\Omega$ and $0 < t < s < d(x, \Omega^c)$, then $$A_f(x; s) - A_f(x; t) = \int h_{s,t,x}\, \rho \,d\lambda$$ where $h \ge 0$ and $\int h \,d\lambda = (s^2 - t^2) / 2(d+2)$. Therefore,
$$|A_f(x; s) - f(x)| = \lim_{t \to 0} |A_f(x; s) - A_f(x; t)| \le {Cs^2 \over 2(d+2)}.$$

Fix $x, y \in \Omega$ with $|y - x| < \frac12 d(x, \Omega^c)$. Suppose $f(x) = 0$. Write
\begin{align*}f(y) &= f(y) - f(x)\\&=[f(y) - A_f(y;s)] + [A_f(y;s) - 2^dA_f(x; 2s)] + 2^d[A_f(x;2s) - f(x)].\end{align*}
The first summand is bounded in absolute value by $Cs^2/2(d+2)$, the second one is nonpositive because $B_s(y) \subseteq B_{2s}(y)$ and $\lambda(B_{2s}) = 2^d \lambda(B_s)$, and the third one is bounded by $2^{d+2} C s^2 / 2(d+2)$. So,
$$f(y) \le Cs^2 {1 + 2^{d+2} \over 2(d+2)}
\le
2^d Cs^2.$$ This is the result. \qed

\begin{theorem}%If $f \ge 0$ is a limit of radial averages, and $\nabla^2(f \,d\lambda) = \rho \,d\lambda$ with $|\rho| \le C$, then $\rho = 0$ a.e. on the set $\{x: f(x) = 0\}$.
Suppose $f \ge 0$ is a limit of radial averages and $\nabla^2(f \,d\lambda) = \rho \,d\lambda$ with $|\rho| \le C$. Then $\rho = 0$ a.e.~on the zero set $\{x \in \Omega: f(x) = 0\}$.\label{zeroonboundary}\end{theorem}

\angry Proof. Let $Z := \{x \in \Omega: f(x) = 0\}$. By the Lebesgue density theorem, there is a set of zero $\lambda$-measure $N$ so that for every point $x \in Z \setminus N$, both of the following equalities hold: \begin{align}\lim_{r \to 0} {\lambda(B_r(x) \cap Z^c) \over \lambda(B_r)} &= 0 \qquad \text{and}\label{zeroonboundaryindicator}\\\lim_{r \to 0} {1 \over \lambda(B_r)} \int_{B_r(x)} \rho \,d\lambda &= \rho(x).\label{zeroonboundarymeasure}\end{align} By the last theorem, $f(y) \le 2^d C |y-x|^2$ for $y$ sufficiently close to $x$, so $$A_f(x; s) \le {\lambda(B_s(x) \cap Z^c) \over \lambda(B_s(x))} \times O(s^2) = o(s^2).$$
This estimate holds for every $x \in Z\setminus N$.

\medskip
For every $x \in Z$ and $s < d(x, \Omega^c)$, we have the inequality $f(x) = 0 \le A_f(x; s)$, so $f$ is subharmonic on average on $Z$. The Laplacian is the signed measure $\rho \,d\lambda$, and by Theorem \ref{superharmonicmeansnegative}, that signed measure must be positive on~$Z$, which means that we must have $\rho \ge 0$ almost everywhere on $Z$.

So it is enough to prove $\rho \le 0$ a.e.~on the zero set. Suppose not. Then there must be at least one point~$x \in Z \setminus N$ with $\rho(x) > 0$. The point is not in $N$, so the limit in equation~\ref{zeroonboundarymeasure} exists. For $s > 0$ sufficiently small, we have
$$\inf_{t \in (0, s)} {1 \over \lambda(B_t)} \int_{B_t(x)} \rho \,d\lambda > {\rho(x) \over 2} \qquad\text{ for } t \in (0, s).$$
Construct $h_{s, t, x}$ as in Theorem \ref{ballcomparison}, and repeat the reasoning in the proof of Lemma \ref{superharmonicmeansnegative} to get the inequality
$$\int h_{s, t, x} \,\rho \,d\lambda \ge {\rho(x) \over 2} \int h_{s, t, x} \,d\lambda.$$
Then by Theorem~\ref{ballcomparison} and the Lemma~\ref{theintegralofh},
$$A_f(x; s) - A_f(x; t) = \int h_{s, t} \,\rho \,d\lambda \ge {\rho(x) \over 2} \int h_{s, t} \,d\lambda = {\rho(x) \over 4(d+2)} (s^2 - t^2).$$
Take the limit of both sides as $t \to 0$ and use the fact that $\lim_{t \to 0} A_f(x;t) = f(x) = 0$ to get the inequality
$$A_f(x; s) \ge {\rho(x) \over 4(d+2)} s^2.$$
This contradicts the estimate $A_f(x; s) = o(s^2)$ for $x \in Z$.

Therefore, $\rho = 0$ a.e.~on the zero set $Z$.\qed

\end{subsection}
\end{section}
\begin{section}{The existence of quadrature domains}
\begin{subsection}{A maximization problem}\label{maximizationproblem}
A weight function is \emph{properly supported} if it is greater than or equal to $1$ on some bounded open set, and $0$ outside that open set.
For example, a sum of indicator functions of bounded open sets is properly supported, but the function $\frac12 \1_{B_1}$ is not properly supported.

From now on, suppose $w$ is a properly supported weight function.
We will prove the existence of a quadrature domain for $w$.
We start by posing a maximization problem, then extract a set from the solution, and finally prove that the set is a quadrature domain in Theorem~\ref{theresone}.

\medskip
\angry Definition. If $\psi$ and $\psi'$ are distributions, then we say that $\psi \le \psi'$ if $\psi' - \psi$ is a positive distribution.

\medskip
We will put a weak continuity condition on the functions in our problem. Let $\emph{RA}(\mathbb R^d)$ be the set of functions $f: \mathbb R^d \to \mathbb R$ which are a limit of radial averages at every point in $\mathbb R^d$. See Section~\ref{lora} for the definition of a ``limit of radial averages.'' These functions are all locally integrable by definition, so $\nabla^2(f d \lambda)$ exists as a distribution. It's a big class of functions: for example, any integrable superharmonic or subharmonic function is in $\emph{RA}(\mathbb R^d)$.

The maximization problem is:

\medskip\noindent
\hskip -.49em\framebox{\hskip .2em \parbox{\hsize}{\vskip .2em
\angry Maximization problem.

Find the largest nonpositive $f \in \emph{RA}(\mathbb R^d)$ with $\nabla^2(f \,d\lambda) \ge (w-1) \,d\lambda$ in the sense of distributions.
\vskip .2em
}\hskip .2em}

\bigskip

It will turn out that there is a function that is pointwise greater than or equal to any other function, so there will really be a largest function.

If we allow the whole class of nonpositive locally integrable functions, we could set any function to zero on a set of measure zero without affecting its distributional Laplacian. That means that no function can be pointwise largest except for $0$, which is typically not a solution.
\end{subsection}
\begin{subsection}{Newtonian potentials}
\label{newtonian}
If $w$ is a bounded, compactly supported weight function, let the {\it Newtonian potential} of $w$ be the convolution of $w$ with the unrestricted Green's function: $$Nw(x) := \int_{\mathbb R^n} G(x, y) w(y) \,dy.$$
Recall that $G(x, y)$ is a function of $x - y$ only, so this is really a convolution.

The convolution of a bounded, compactly supported function with a locally integrable function is continuous.
Green's function is locally integrable, so $Nw$ is continuous on all of $\mathbb R^d$. It has the first derivative $${\partial \over \partial x_i} Nw(x) = \int_{\mathbb R^n} {\partial G(x,y) \over \partial x_i} w(y) \,dy = -{1 \over C_d} \int_{\mathbb R^n} {x_i - y_i \over |x-y|^d} w(y) \,dy$$ which is continuous on $\mathbb R^d$ for the same reason.
To see that this is really the derivative, one can integrate, use Fubini's theorem, and check that the result is $Nw$ plus a function independent of $x_i$.

We have already observed in Corollary~\ref{dist} that $-\nabla^2 (N\varphi) = \varphi$ in the sense of distributions, so $-\nabla^2 (N w d\lambda) = w d\lambda$.
If $w$ is bounded and {\it constant} outside a compact set, we can still define $Nw$: let $c$ be the constant, and set $Nw := N(w - c) + {c \over 2d}|x|^2$. Again $-\nabla^2 (Nw \,d\lambda) = (w-c) d\lambda + c d\lambda = w d\lambda$.

Note that a continuous function is necessarily a limit of radial averages, so if $f$ is a limit of radial averages, then so is $f + Nw$.

\end{subsection}
\begin{subsection}{Maximization over subharmonic functions}
\def\tle{\mathrel{\trianglelefteq}}

\begin{theorem}\label{minimizationequivalent}The maximization problem is equivalent to:

{\narrower

\medskip
\noindent
Find the largest nonpositive function $f$ on $\mathbb R^d$ with the property that the sum $f + N(w-1)$ is subharmonic everywhere in $\mathbb R^d$.

}\end{theorem}

\angry Proof.
We will show that $f$ is a limit of radial averages with $\nabla^2 (f \,d\lambda) \ge (w-1)\,d\lambda$ if and only if $f + N(w-1)$ is subharmonic, so the set of admissible functions is the same for both problems.

\medskip

Suppose $f$ is a limit of radial averages and $\nabla^2 (f\,d\lambda) \ge (w-1)\,d\lambda$. Then $$\nabla^2[(f + N(w-1)) \,d\lambda] = \nabla^2 (f \,d\lambda) - (w-1) \,d\lambda \ge 0.$$ By Lemma~$\ref{negativemeanssuperharmonic}$, there is a subharmonic $\bar f$ which is equal almost everywhere to $f + N(w-1)$, but in fact this equality holds everywhere, because both sides are a limit of radial averages. Therefore, $f+N(w-1)$ is subharmonic, and $f$ is an admissible function for the second problem.

\medskip
On the other hand, suppose $f + N(w-1)$ is subharmonic. By Theorem~$\ref{superharmonic}$, $\nabla^2[(f + N(w-1)) \,d\lambda] = \nabla^2(f d\lambda) - (w-1) d\lambda \ge 0$, so $\nabla^2(f d\lambda) \ge (w-1) d\lambda$. Finally, $f$ is the difference of a subharmonic function and a continuous function, $f = (f+N(w-1)) - N(w-1)$, so it is a limit of radial averages, and it is admissible for the first problem.

The admissible functions for both problems are the same, so they are equivalent.\qed

\medskip
We can now use the fundamental convergence theorem for subharmonic functions to find a minimum.

\begin{theorem}[Fundamental convergence theorem]\

Let $\Gamma$ be a family of subharmonic functions defined on an open subset of $\mathbb R^d$ and locally uniformly bounded above.
Let $u(x)$ be the pointwise supremum of all the functions in $\Gamma$. Let $u_+(x) = \max\{u(x), \limsup_{y \to x} u(y)\}$.

Then~$u_+ = u$ almost everywhere, and $u_+$ is subharmonic.
\end{theorem}

\angry Proof. See for example Section 1.III.3 of Doob \cite{doob} although this is stated in terms of superharmonic functions.

\begin{corollary} Let $\pi$ be a measurable, bounded function on $\mathbb R^d$ that is constant outside a compact set. Then there is a largest $f \le 0$ with the property that $f + N\pi$ is subharmonic. \label{smallestsuperharmonic}
\end{corollary}

\angry Proof.
Let $\Gamma = \{u: u\text{ is subharmonic},\, u \le N\pi\}.$
The functions in this class are uniformly bounded above on any compact set $K$ by $\max_K N\pi$.

Apply the fundamental convergence theorem to $\Gamma$ to get a subharmonic function $u_+$ greater than or equal to every function in $\Gamma$. Then $u_+ \in \Gamma$:
\begin{align*}u_+(x) &= \max\{u(x), \limsup_{y \to x} u(y)\} \\&\ge \max\{N\pi(x), \limsup_{y \to x} N\pi(y)\} \\&\ge N\pi(x)\end{align*} by continuity of $N\pi(x)$, so it satisfies the inequality and is subharmonic.

Set $f := u_+ - N\pi$. Then $f \le 0$ and $f + N\pi$ is subharmonic. If $g$ is any other function with $g \le 0$ and $g + N\pi$ subharmonic, then the sum $g + N\pi$ is in the class $\Gamma$, so $u_+ \ge g + N\pi$ and $f \ge g$. \qed

\begin{corollary}
\label{minimizable}
There is a largest function $f \le 0$ that is a limit of radial averages and satisfies $\nabla^2(f \,d\lambda) \ge (w-1) \,d\lambda$.
\end{corollary}

\angry Proof. Combine Corollary~\ref{smallestsuperharmonic} with Theorem~\ref{minimizationequivalent}.\qed

\medskip
In the next section, we will characterize the Laplacian of the minimal function, and discover that there is a quadrature domain hiding inside it.

\end{subsection}
\begin{subsection}{Finding the Laplacian}

In the rest of this section, we fix a properly supported weight function $w$, and we define an associated open set $A$ and a locally finite measure $\mu$.

Let $f \ge 0$ be the solution of the maximization problem in Section~\ref{maximizationproblem}. Let $u := f + N(w-1)$, which is a subharmonic function.

Observe that $f = u - N(w-1)$ is the difference of a subharmonic function and a continuous function, so it's upper semicontinuous. This means that the sets $\{f < c\}$ are open for $c \in \mathbb R$.

\medskip
\angry Definition. Let $A$ be the open set $\{x:f(x) < 0\}$.

\medskip We will prove later that $A \cup \{w \ge 1\}$ is the quadrature domain we are looking for. See Corollary \ref{nnlezero} and Theorem \ref{theresone}.

A subharmonic function has a positive distributional Laplacian, which is a locally finite measure by Corollary \ref{positivemeansmeasure}. We give this measure a name.

\medskip
\angry Definition. Let $\mu$ be the locally finite measure with $\nabla^2 (u \,d\lambda) = d\mu$.

\medskip
We now study the geometry of $\mu$ and $A$. We first want to prove the basic fact that $\mu$ is supported on $A^c$.

\begin{lemma}
$\mu(A) = 0$.%
\label{harmoniconpositive}
\end{lemma}

\angry Proof.
Let $x \in A$, so that we have $f(x) < 0$ and $u(x) < N[w-1](x)$.
Choose a small radius $r>0$ with $\overline{B_r(x)} \subseteq \Omega$ and $$\max_{y \in \overline{B_r(x)}} u(y) < \min_{y \in \overline{B_r(x)}} N[w-1](y).$$ This is possible because $u$ is upper semicontinuous.

We can use the Poisson kernel to construct a function $u'$ that is equal to $u$ outside $B_r(x)$, and is harmonic on $B_r(x)$. The resulting function will still be subharmonic, and $u' \ge u$. See e.g.~I.II.6~of Doob~\cite{doob}.

Then $u' \le N(w-1)$, because $u'(y) \le \max_{B_r} u \le \min_{B_r} N(w-1)$ inside the ball and $u' = u \le N(w-1)$ outside the ball. But $u$ is the largest subharmonic function that is less than or equal to $N(w-1)$, so $u' = u$.

It follows that the original function $u$ was harmonic on $B_r(x)$. If a function is harmonic on an open set, then its Laplacian is zero on that set, so $\nabla^2(u \,d\lambda)[h] = 0$ for any function $h \in C^\infty_c(B_r(x))$. Therefore $\mu(B_r(x)) = 0$. We can cover $A$ by countably many balls $B_r(x)$, so $\mu(A) = 0$.\qed

\medskip

We now bring in the theory from Section~\ref{pos} to find an explicit formula for the distributional Laplacian $\nabla^2(f \,d\lambda)$ in terms of $w$ and $A$.
First, we'll use Lemma~\ref{superharmonicmeansnegative} together with the Lebesgue density theorem to get bounds on the restriction of $\mu$ to $\partial A$.

\begin{theorem}If $E \subseteq A^c$ is measurable, then $\mu(E) \le \lambda(E)$.\label{absolutecontinuity}\end{theorem}

\angry Proof.
We know $f$ is a limit of radial averages. If $x$ is any point in $A^c$, then $f(x) = 0$, so $f(x) \ge A_r(f; x)$ just because the function is nonpositive.
Therefore, $-f$ is subharmonic on average on $A^c$.

Let $\nu$ be the signed measure $\nu(E) := \mu(E) + \int_E (w-1) \,d\lambda$. Then $d\mu = \nabla^2(u \,d\lambda) = \nabla^2(f \,d\lambda) - (w-1) \,d\lambda$ by definition and Corollary~\ref{dist}, so
\begin{align*}
\nabla^2(f d\lambda) &=
    \nabla^2(u \,d\lambda) + (w-1) \,d\lambda
    \\&= d\mu + (w-1) \,d\lambda
    \\&= d\nu.
\end{align*}
Of course, that means that the distributional Laplacian of $-f$ is also a signed measure, $\nabla^2(-f \,d\lambda) = -d\nu$.
This allows us to use
Lemma~\ref{superharmonicmeansnegative} on $-f$, and we get $-\nu(E) \ge 0$ for measurable $E \subseteq A^c$.

The weight function is nonnegative by definition, so $$\mu(E) = \nu(E) - \int_E (w-1) \,d\lambda \le \int_E 1 \,d\lambda = \lambda(E).$$
This proves the result.
\qed

\begin{corollary} $\nabla^2(f \,d\lambda) = (w-1) \1_A \,d\lambda$.\label{findthelaplacian}
\end{corollary}
The theorem implies that $\mu$ is absolutely continuous with respect to $\lambda$, so there is a Radon-Nikodym derivative $d\mu / d\lambda$.

Let $\rho = d\mu / d\lambda + w - 1$, so that $\rho \,d\lambda = d\nu = \nabla^2(f \,d\lambda)$. This function is bounded, because $-1 \le \rho \le \max w$, so we are in the setting of Theorem \ref{zeroonboundary}, which tells us that $\rho = 0$ a.e.~on the set where $f$ is zero, which is $A^c$.

Theorem \ref{harmoniconpositive} says that the measure $\mu$ is zero on $A$, so $\rho = w-1$ on $A$ a.e. Therefore, $\nabla^2(f \,d\lambda) = \rho \,d\lambda = (w-1) \1_A \,d\lambda$.
\qed

\end{subsection}
\begin{subsection}{Quadrature domain inequality for Green's functions}

In this section, we'll prove that $Q := A \cup \{w \ge 1\}$ satisfies the quadrature domain inequality as long as $h$ is one of Green's functions.

\begin{theorem}If $f$ solves the maximization problem, then $f = N\big[(1-w)\1_A\big]$.\end{theorem}

\angry Proof. Let $\rho := (w - 1) \1_A$.
By Corollary~\ref{dist}, $-\nabla^2(N \rho \,d\lambda) = \rho \,d\lambda$, so $$\nabla^2[(f - N\rho) \,d\lambda] = 0.$$
Then $f - N\rho$ is a harmonic function.\footnote{As before, we use Corollary \ref{negativemeanssuperharmonic} to get a subharmonic function $\bar f$ that is equal a.e., and then argue that $f - N\rho$ is a limit of radial averages, so it has to be exactly equal to $\bar f$.}

We want to show that it's zero. Our strategy will be to prove that it goes to zero at infinity and use the maximum principle.

We start with the remark that the set $A = \{f < 0\}$ is always bounded. To prove this, one just has to exhibit a function $h$ that is admissible in the maximization problem and is compactly supported. This is straightforward but somewhat tedious and we postpone it to the next section.

If $d \ge 3$, then the unrestricted Green's function is $|x-y|^{2-d} / (d-2)C_d$. Because $A$ is bounded, this
goes to zero as $x \to \infty$ uniformly for $y \in A$.
So $$N\rho(x) = \int G(x, y) \rho(y) \,dy$$ converges uniformly to zero as $x \to \infty$. Of course, $f$ is compactly supported, so $f(x) \to 0$ as $x \to \infty$.
Now we can apply the maximum principle: if $|x| \le r$,$$\min_{z \in \partial B_r} f(z) - N\rho(z) \le f(x) - N\rho(x) \le \max_{z \in \partial B_r} f(z) - N\rho(z),$$ and taking limits as $r \to \infty$, we get $f(x) - N\rho(x) = 0$ for $x \in \mathbb R^d$.

What about when the dimension is smaller? In $d = 2$, we'll use a similar strategy, but we require more information than before. Let $h$ be a smooth, compactly supported function which is $1$ on $A$.
Then
$$\int h \rho \,d\lambda = (\rho \,d\lambda)[h] = (f \,d\lambda)[\nabla^2 h] = \int f \nabla^2 h \,d\lambda = 0,$$
because $\nabla^2 h = 0$ on $A$ and $f = 0$ outside $A$. But $\rho = h \rho$, so $\int \rho \,d\lambda = 0$. This means in particular that the integral against any constant is zero.

Then we again have
\begin{align*}N\rho(x) &= \int G(x, y) \rho(y) \,dy \\&= \int (G(x, y) - G(x, 0)) \rho(y) \,dy \\&= \int {1 \over \pi} \log {|x| \over |x-y|} \rho(y) \,dy \to 0.\end{align*}

In one dimension, we also want to know that $\int y \rho \,dy = 0$, which can be proved in the same way using $h \in C^\infty_c$ with $h(y) = y$ for $y \in A$. Then $$N\rho(x) = \int G(x, y) \varphi(y) \,dy = \int (y-x) \rho(y) \,dy$$
for $x \ge \max A$, and this is zero because $\int \rho \,dy = \int y \rho \,dy = 0$. Similarly, $N\rho(x) = 0$ when $x \le \min A$. In each case, we conclude that $f = N\rho$.\qed

\medskip
We restate this in terms of a set $Q$ which will turn out to be our quadrature domain.

\begin{corollary}Let $Q := A \cup \{w \ge 1\}$. Then $f = N[\1_Q - w]$.\label{qco}
\end{corollary}

\angry Proof.
By the proof of Corollary \ref{findthelaplacian}, $\nu$ is zero on $A^c$. Let $E$ be any measurable subset of $A^c$. Then $\mu(E) = \nu(E) -\int_E (w-1) \,d\lambda = -\int_E (w-1) \,d\lambda$. For this to be nonnegative, we must have $w \le 1$ almost everywhere in $A^c$.

Part of the definition of a good weight function is that $w = 0$ or $w \ge 1$ at every point, which means that $\1_Q = w$ a.e.~on $A^c$. Therefore, $$\1_Q - w = (\1_Q - w)\1_A = (1-w)\1_A \qquad\text{a.e.}$$ It follows immediately that $f = N[(1-w)\1_A] = N[\1_Q - w]$.
\qed

\medskip
We now prove that $Q$ satisfies the quadrature domain inequality if $h$ is one of the Green's functions.

\begin{corollary}\label{nnlezero}
If $\varphi = \1_Q - w$, then $N\varphi \le 0$ on $Q$ and $N\varphi = 0$ on $Q^c$.
If $h$ is subharmonic on $Q$ and $h = \pm G_x$ for $x \in \mathbb R^d$, then $$\int hw \,dy \le \int_Q h \,dy.$$
\end{corollary}

The first statement is obvious: we have proven that $f = N\varphi$, and $f$ is nonnegative everywhere by definition, so $f(x) \le 0$ on $Q$. If $x \in Q^c$, then we are outside the set $A = \{f < 0\}$, so $f(x) = 0$.

The second statement follows almost immediately.
By the definition of the Newtonian potential, $$\int_Q G_x \,dy - \int G_x w \,dy = N[\1_Q - w](x) = f(x).$$

The functions of the form $\pm G_x$ that are subharmonic on $Q$ are $-G_x$ for all points $x \in \mathbb R^d$, and $G_x$ for $x \notin Q$. If we set $h = -G_x$, then $$\int_Q h \,dy - \int hw \,dy = -f(x) \ge 0$$ for any point $x \in \mathbb R^d$. In the same way, if $h = G_x$, then $\int_Q h \,dy - \int hw \,dy = f(x)$, which is zero for any point $x \in Q^c$. That proves the corollary.\qed

\medskip

In the proof above, we used the fact that $A$ is always bounded, and we will prove that in the next section.

After that, we'll start working on Theorem~\ref{sakaislemma}, which states that, if the quadrature domain inequality holds for subharmonic functions of the form $\pm G_x$, then it holds for any integrable subharmonic function.
\begin{subsubsection}{The set $A$ is always bounded}\label{thesetbounded}If $w$ is any properly supported weight function, and $f$ is the solution of the maximization problem, then $A = \{f < 0\}$ is always bounded.

We will prove this by finding a function $\varphi \le 1-w$ where $N{}\varphi$ is compactly supported. First we compute the Newtonian potential of the ball.\begin{lemma} $$N\1_{B_R} = \begin{cases}c_1-|x|^2/2d&|x| \le R\\G(x, 0) \lambda(B_R)&|x| \ge R\end{cases},$$ where $c_1$ is the constant that makes this continuous. \end{lemma}

\angry Proof. If $|x| \ge R$, then the function $y \mapsto G(x, y)$ is harmonic on $B_R$, so $\int G(x, y) \1_{B_R} \,dy$ is equal to $G(x,0)$ times the measure of $B_R$.

 Let $h(x) := N\1_{B_R} + |x|^2/2d$. This function is continuous. The distributional Laplacian of $h$ is zero inside the ball, because $\nabla^2 N \1_{B_R} = -\1_{B_R}$ and $\nabla^2 |x|^2 = 2d$.
By Lemma~\ref{negativemeanssuperharmonic} applied to both $h$ and $-h$, it is both subharmonic and superharmonic, so $h$ is harmonic inside the ball.

The Newtonian potential of the ball is radially symmetric, and so is $|x|^2 / 2d$, so the value of $h$ on the sphere $|x| = R$ is a constant. By the maximum principle, $h$ is constant inside $B_R$, so $N\1_{B_R} = c_1 - |x|^2/2d$ in $B_R$.
\qed

\medskip
Note that the radial derivative of this potential is $${\partial \over \partial r} N\1_{B_R} = \begin{cases}-r/d&|x| \le R\\R^d/dr^{d-1}&|x| \ge R\end{cases}.$$

\begin{corollary}\label{bou}If $w$ is a properly supported weight function, then $A$ is always bounded.\end{corollary}

\angry Proof. Let $c = 1 \vee \max w$. A properly supported weight function is zero outside a bounded set. Let $R$ be the radius of that set, and let $R' = c^{1/d} R$.

Set $\varphi = \1_{B_{R'}} - c\1_{B_r}$. Note $\varphi \le 1-w$, so $$\nabla^2 (N\varphi \,d\lambda) = -\varphi \,d\lambda \ge (w-1) \,d\lambda.$$ Therefore, $N\varphi$ will be admissible in the maximization problem as long as $N\varphi \le 0$.
The lemma tells us that
$$N\varphi = \begin{cases}(c-1)|x|^2/2d - c_1&|x| \le R\\- cG(x,0)\lambda(B_R)- |x|^2/2d - c_2&R \le |x| \le R'\\0&|x| \ge R'.\end{cases}$$ Here $c_1, c_2$ are constants that make this continuous. Note that the function is zero outside the ball $B_{R'}$.
We computed the radial derivative of $N\1_{B_R}$ above, and if we plug it in, we will get
$${\partial \over \partial r} N\varphi = \begin{cases}(c-1)r/d&r \le R\\r(cR^d/r^d - 1)/d&R \le r \le c^{1/d}R\\0&r \ge c^{1/d}R\end{cases}$$ where $r = |x|$. This derivative is at least zero, because $R > r$, so $N\varphi \ge 0$.

Therefore, $N\varphi$ is admissible in the maximization problem that we solved to get $f$, so $0 \ge f \ge N\varphi$ and in particular $f = 0$ whenever $N\varphi = 0$. But $N\varphi = 0$ outside $B_{R'}$, so $A$ must be bounded.\qed

\end{subsubsection}
\end{subsection}
\begin{subsection}{Quadrature domain inequality generally}\label{allintegrablesubharmonicfunctions}
Here is the grand climax of this chapter, Sakai's Lemma~5.1 \cite{sakaiobstacle}. We present the proof and some minor results that are used.

\begin{theorem}[Sakai's Lemma 5.1] Let $Q$ be an open bounded set. If there is a function $\varphi \in L^\infty(Q)$ with $N\varphi \le 0$ on $Q$ and $N\varphi = 0$ on $Q^c$, then
$$\int s\varphi \,dy \ge 0$$
for every integrable subharmonic function $s$ on $Q$.

\label{sakaislemma}\end{theorem}

\medskip
\angry Proof. The basic idea is approximation, but it's delicate and relies on a tight estimate of $N\varphi$ near the boundary of $Q$.

Without loss of generality, we can assume that $|\varphi| \le 1$ everywhere.
Let $s$ be an integrable subharmonic function on $Q$.
Let $s_{n} := s * \psi_n$ be the approximations defined in Theorem~\ref{superharmonic}. As before, each function $s_n$ is defined on $\{d(x, Q^c) > 1/n\}$, and it's smooth and subharmonic.

Let $h_j$ be the sequence of functions from Theorem \ref{lemmafour}. These are smooth functions that are compactly supported in $Q$ with $0 \le h_j \le 1$ that converge pointwise to~$1$, and $h_j \equiv 1$ on the set $\{d(x, Q^c) > 1/j\}$. They have some other properties that we will need.

By the dominated convergence theorem,
$$\lim_{j \to \infty} \int s h_j \varphi \,dy = \int s \varphi \,dy.$$

Each $h_j$ is supported on a compact subset $Q_j$ of $Q$, so $\int sh_j \varphi \,dy$ is the limit of $\int s_n h_j \varphi \,dy$. Therefore,
$$\lim_{j \to \infty} \left(\lim_{n \to \infty} \int s_n h_j \varphi \,dy\right) = \int s \varphi \,dy.$$
If $j$ is fixed, then for large enough $n$, $s_n h_j$ is smooth and compactly supported in $Q$. By Lemma~\ref{test} and the fact that Green's function is symmetric,
$$\int s_n h_j \varphi \,dy = \int N (-\nabla^2) (s_n h_j) \varphi \,dy = \int (-\nabla^2) (s_n h_j) N \varphi \,dy.$$
A smooth subharmonic function has a positive Laplacian, so $\nabla^2 s_n \ge 0$, and by our assumption, $N \varphi \le 0$.
So $\nabla^2(s_n) h_j N\varphi$ is less than or equal to zero. That gives us the inequality
\begin{align*}
(-\nabla^2)(s_nh_j) N \varphi &\ge (-\nabla^2)(s_n h_j) N \varphi + \nabla^2(s_n) h_j N\varphi\\&=-2(\nabla s_n) \cdot (\nabla h_j) N \varphi - s_n \nabla^2 (h_j) N \varphi\\&=-2 \nabla \cdot (s_n \nabla(h_j) N\varphi) + 2s_n (\nabla h_j \cdot \nabla N \varphi) + s_n \nabla^2 (h_j) N\varphi. \end{align*}

Divergences go away when we integrate, because all these functions are smooth and compactly supported in $Q$, and we get a lower bound
\begin{align}
    \int s_n h_j \varphi \,dy
        &= \int(-\nabla^2)(s_nh_j)N\varphi\,dy\notag
        \\&\ge\int2s_n(\nabla{h_j}\cdot\nabla{N}\varphi)
            +s_n\nabla^2(h_j)N\varphi\,dy.\label{twointegral}
\end{align}

We have managed to get the dependence on $n$ outside of the derivatives. Our goal now is to prove a lower bound of the form $-C/j$ on (\ref{twointegral}) and transfer it to $\int s \varphi \,dy = \lim_j \lim_n \int s_n h_j \varphi \,dy$.

The test functions from Theorem \ref{lemmafour} are constant in the set $\{d(x,Q^c) > 1/j\}$, so $\nabla h_j = 0$ and $\nabla^2 h_j = 0$ in that set. So to bound the integrals in (\ref{twointegral}), we need to study $N\varphi, \nabla N \varphi, \nabla h_j, \nabla^2 h_j$ near the boundary of $Q$.

\medskip
Recall that $N \varphi$ has a continuous first derivative. For $y \in Q$, let $\delta(y) = d(y, Q^c)$. Our assumption says that $N \varphi = 0$ on $Q^c$, and the function isn't positive anywhere, so $\nabla N \varphi$ must be zero on $Q^c$.

By Lemma~\ref{continuityofnewtonian} below, there are constants $\delta_0 > 0$ and $C > 0$ with $$\left|\nabla N \varphi(y) - \nabla N \varphi(x)\right| \le C |y-x| \log {1 \over |y-x|} $$ for $y, x \in \mathbb R^d$ with $|y-x| < \delta_0$. We have seen that $\nabla N \varphi = 0$ outside $Q$, so the gradient is small near the boundary: $|\nabla N \varphi(y)| \le C \delta \log 1/\delta$ for $\delta < \delta_0$.

The bound on the gradient implies a bound on $N\varphi$:
\begin{align*}|N \varphi(y)|
    &= |N\varphi(y) - N\varphi(x)|
    \\&\le \int_0^\delta Ct \log 1/t \,dt
    \\&\le \int_0^\delta C(2t \log 1/t - t) \,dt
    \\&= C \delta^2 \log 1/\delta.
\end{align*}
In the third line, we assume that $\delta_0 < 1/e$, so $\log 1/t > 1$.

We chose test functions $h_j$ from Lemma \ref{lemmafour}, so we can use the bounds in the lemma. The magnitude of the gradient is $$|\nabla h_j| \le {1 \over j \delta \log 1/\delta}$$ and $|\nabla^2 h_j| \le |\partial^2h_j/\partial{y}_1^2| + \cdots + |\partial^2h_j/\partial{y}_d^2| \le 2d/(j \delta^2 \log 1/\delta)$. These blow up as $y$ approaches the boundary of $Q$, but they grow slowly enough that the bounds on $N\varphi$ and $\nabla N \varphi$ cancel them out.

We have $|\nabla h_j| |\nabla N \varphi| \le C / j$, which gives a bound on the first part of (\ref{twointegral}):\begin{align*}\left|\int 2s_n (\nabla h_j \cdot \nabla N \varphi) \,dy\right| &\le {2C \over j} \int_{Q_j} |s_n| \,dy\end{align*}

For the second integral, we have the bound $|\nabla^2 h_j| |N \varphi| \le 2Cd/j$, so $$\left|\int s_n \nabla^2 (h_j) N \varphi \,dy\right| \le{2Cd \over j} \int_{Q_j} |s_n| \,dy.$$

Set $C' = 2(1+d)C$. Then we can put these bounds together to get $$\int s_n h_j \varphi \,dy \ge -{C' \over j} \int_{Q_j} |s_n| \,dy.$$

Again, $|s|$ is integrable, so if we take the limit on both sides and use bounded convergence, we get $\int s h_j \varphi \,dy \ge -(C'/j) \int |s| \,dy$.

Finally, we can take the limit as $j \to \infty$, and this shows that $\int s\varphi \,dy \ge 0$: $$\int s \varphi \,dy \ge \lim_{j \to \infty} \int s h_j \varphi \,dy \ge \lim_{j \to \infty} -{C' \over j} = 0.$$ This is what we were trying to prove.\qed

\medskip
The corollary is that quadrature domains exist.
\begin{theorem}Every properly supported weight function has a quadrature domain.\label{theresone}\end{theorem}

\angry Proof. We have a recipe for the quadrature domain: we solve the maximization problem in Section \ref{maximizationproblem}, and then we set $Q := \{f < 0\} \cup \{w \ge 1\}$ as in Corollary \ref{qco}. This is an open set, and it's bounded, by Lemma~\ref{bou}.

By Corollary \ref{nnlezero}, $f = N(\1_Q - w)$ is nonpositive everywhere and zero outside $Q$, so we can use Theorem \ref{sakaislemma} with $\varphi = \1_Q - w$. This tells us that for any integrable subharmonic function $s$ on $Q$,
$\int s(\1_Q - w) \,dy \ge 0$.

In other words, if  $s$ is integrable and subharmonic on $Q$, then
$$\int sw \,dy \le \int_Q s\,dy.$$
That's the definition of a quadrature domain.\qed

\end{subsection}
\begin{subsection}{Denouement 1: log-Lipschitz continuity}%We've used two lemmas that we must prove:
We owe two lemmata that we must prove. First, a lemma which says that the first derivative of the Newtonian potential is a little worse than Lipschitz.

\begin{lemma}[G\"unther, \cite{guntherpotentialtheory}, \S13]\label{continuityofnewtonian}Suppose $\varphi$ is bounded and measurable, and zero outside a bounded set $E$. If $y, y' \in \mathbb R^d$ and $|y - y'| = \varepsilon$, then $$\left|{\partial N\varphi \over \partial y_i}(y) - {\partial N\varphi \over \partial y_i}(y')\right| = O\!\left(\varepsilon \log {1 \over \varepsilon}\right).$$
The constant in the $O$-notation depends only on $\max |\varphi|$ \!and $\mathop{\rm diam} E$.
\end{lemma}

\angry Proof. Write both terms as derivatives of integrals, move them both under the same integral sign, and move the derivative inside the integral, to get
\begin{align*}\left|{\partial N\varphi \over \partial y_i}(y) - {\partial N\varphi \over \partial y_i}(y')\right|
= \left|\int_{\Omega} \left[{\partial G \over \partial y_i}(x, y) - {\partial G \over \partial y_i}(x, y')\right] \varphi(x) \,dx\right|.
\end{align*}
This is justified because $\partial G / \partial y_i$ is locally integrable and $\varphi$ is bounded.

Let $A = \{x \in \Omega: |x - y| < 2 \varepsilon\}$, and break the integral on the last line up into $\int_A$ and $\int_{E \setminus A}$.
The first derivatives of $G(x, y)$ are~$\smash{O(||x-y||^{1-d})}$, so~$$\int_A \left[{\partial G \over \partial y_i}(x, y) - {\partial G \over \partial y_i}(x, y')\right] \varphi(x) \,dx = O\left(\int_{B_{2 \varepsilon}(y)} ||x-y||^{1-d} \,dx\right) = O(\varepsilon).$$
\vskip .3em
For the part outside of $A$, we estimate the integrand with derivatives.
By the mean value theorem,
there is a point $y''$ on the line segment between~$y$ and~$y'$ with ${\partial G \over \partial y_i}(x, y) - {\partial G \over \partial y_i}(x, y') = \smash{(y - y') \cdot \nabla {\partial G \over \partial y_i} (x, y'')}$.
The second derivatives of $G(x; y)$ are $O(||x-y||^{-d})$, so that dot product is at most
$$|y - y'| \times \left|\nabla {\partial G \over \partial y_i}(x, y'')\right| =
O(\varepsilon ||x-y''||^{-d}).$$

When $x$ is not in $A$, $||x-y|| \ge 2\varepsilon$, so \begin{align*}||x - y''|| &\ge ||x - y|| - ||y - y''|| \\&\ge ||x - y|| - ||y-y'|| \ge {\textstyle\frac12} ||x-y||.\end{align*}
Therefore $||x - y''||^{-d} \le 2^d ||x - y||^d$ and $O(||x - y''||^{-d}) = O(||x - y||^{-d})$, and the dot product is bounded by $|(y-y') \cdot \nabla {\partial G \over \partial y_i}(x, y'')| = O(\varepsilon ||x-y||^{-d})$, where the constants in the $O$-notation depend only on $\max |\varphi|$.

If $E \setminus A$ is empty, the integral over it is zero. Otherwise, let $\overline{B_{r_1}} \setminus B_{r_0}$ be the minimal closed annulus containing $E \setminus A$, i.e. $r_0 := \inf\{|x-y|: x \in E \setminus A\}$ and $r_1 := \sup\{|x-y|: x \in E \setminus A\}$, with $r_0 \ge 2\varepsilon$ and $r_1 \le \mathop{\rm diam} E$.

We can estimate the integral over $E\setminus A$ by
\begin{align*}
\int_{E \setminus A} \left|(y - y') \cdot \nabla {\partial G \over \partial y_i}(x, y'')\right| \, \varphi(x) \,dx &= \int_{E \setminus A} O(\varepsilon ||x - y||^{-d}) \,dx
\\&\le O\!\left(\int_{r_0}^{r_1} \varepsilon r^{-d} \,C_d r^{d-1} \,dr\right) \\&= O\!\left(\varepsilon \log {r_1 \over r_0}\right).\end{align*}
Again, the constants in the $O$-notation depend only on $\max |\varphi|$.

Then $r_1 / r_0 \le \mathop{\text{diam}} E/2\varepsilon$, so the estimate is $O(\varepsilon \log r_1 / r_0) = O(\log 1/\varepsilon)$, where the constants in the notation now depend on $\max |\varphi|$ and $\mathop{\text{diam}} E$.
Combine the estimates for $A$ and $E \setminus A$ to get the result:
\begin{align*}\int_E \left[{\partial G \over \partial y_i}(x, y) - {\partial G \over \partial y_i}(x, y')\right] \varphi(x) \,dx &= \int_A (\cdots\hskip -.006em) \varphi(x) \,dx + \int_{E \setminus A} (\cdots\hskip -.006em) \varphi(x) \,dx \\&=
O(\varepsilon) + O\!\left(\varepsilon \log {1 \over \varepsilon}\right)
\\&=O\!\left(\varepsilon \log {1 \over \varepsilon}\right).\end{align*}
This is our bound on the derivatives of the Newtonian potential.
\qed

\medskip
\begin{corollary}The same holds if $\varphi$ is constant outside a bounded set.\end{corollary}

\angry Proof. If $\varphi(x) = c$ as $x \to \infty$ then $N\varphi = N(\varphi -c ) + c|x|^2/2d$. The first term has $O(\varepsilon \log 1/\varepsilon)$ modulus of continuity, and $|x|^2$ is Lipschitz.\qed

\end{subsection}
\begin{subsection}{Denouement 2: a smooth approximation of the distance function}

Let $\delta(x) = d(x,Q^c)$. The second lemma that we owe gives us smooth test functions $h_j$ that are constant on $\{\delta(x) > 1/j\}$ and which have good pointwise bounds on the derivatives. For example,
$$|\nabla h_j(x)| \le {1 \over j \delta \log(1/\delta)}.$$

Compare this to the naive choice $\hat h_j = 1 \wedge jd(x,Q^c)$ where $\max_j |\nabla \hat h_j(x)|$ is of order $1/\delta$ (ignoring smoothness issues).
We are able to reduce the pointwise bound by a factor of $j \log(1/\delta)$ compared to this naive function.

Let $\xi(t) := 1/(t \log 1/t)$. Then $\xi$ is decreasing on $(0, 1/e)$, and for any small $\varepsilon$, the integrals of $\xi$ and $|\xi'|$ on $(0, \varepsilon)$ are $+\infty$.
\begin{lemma}\label{somerealvariablenonsense} For any $m > 0$ there is a function $\eta_m \in C^\infty_c(0, 1/m)$ with $\int_0^\infty \eta_m \,dt = 1$ where $0 \le \eta_m \le \xi / m$ and $|\eta'_m| \le \xi/mt$.\end{lemma}

\angry Proof. Let $m > e$ without loss of generality, so $\xi' \le 0$.
Choose a sequence of nonnegative smooth functions $f_n$ on $(0, 1/m)$ that increase to $-\xi' / m \ge 0$.

Let $F_n(t) := \int_t^{1/m} f_n(\tau) \,d\tau$. Then $F_n(t) \le \xi(t) / m$, and the integral of $\xi(t)$ on $(0, 1/m)$ is $+\infty$, so $I_n := \int_0^{1/m} F_n(t) \,dt \to \infty$ as $n \to \infty$.

Choose $n$ with $I_n \ge 1$. Let $\eta_m(t) := F_n(t) / I_n$. Then this function has the desired properties: it's smooth and compactly supported on $(0, 1/m)$, the integral over that interval is $I_n/I_n = 1$, the function is bounded above by $F_n(t) \le \xi(t)/m$, and the derivative of $\xi$ is negative and bounded by
$$-\xi'(t) = \left(\log {1 \over t} - 1\right) \xi^2(t) \le \left(\log {1 \over t}\right)\xi^2(t) = {\xi(t) \over t},$$
so $|\eta_m'| = |F_n'(t) / I_n| \le -\xi'/ m \le \xi / mt$.\qed

\begin{lemma}[Stein \S VI.2 Theorem 2 \cite{steinsingular}]\ \label{partitionofunitygivesdistancefunction}
Let $Q$ be an open set. There exists a smooth function $\Delta_Q(x)$ on $Q$ with $d(x, Q^c) \le \Delta_Q(x) \le C d(x, Q^c)$ for some positive constant $C > 0$ depending only on the dimension, and
$$\left|\nabla \Delta_Q\right| \le C \qquad \left|{\partial^2 \over \partial x_a \partial x_b} \Delta_Q\right| \le {C \over d(x, Q^c)}.$$
\end{lemma}

\angry Sketch of proof. Let the side length of a cube $\omega$ be denoted by $\mathop{\rm side}(\omega)$. The set $Q$ can be written as a union of closed cubes with disjoint interiors in such a way that the side length of each cube $\omega$ is a power of two, and $$d(\omega, Q^c) \le 8 \mathop{\rm side}(\omega) \le 4d(\omega, Q^c).$$

Scale up each cube around its center by a factor of $1 + 1/\sqrt{d}$. All the scaled cubes are still contained in $Q$.

Pick a smooth function $h \ge 0$ that is compactly supported in the open cube with side $1 + 1/\sqrt{d}$ centred at the origin and so that $h \equiv 1$ inside the unit cube centred at the origin. If $y_\omega$ is the centre of $\omega$, let $h_\omega(x) = h((x - y_\omega) / \mathop{\rm side}(\omega))$. Then $h_\omega \equiv 1$ on $\omega$, and its support is contained inside the scaled cube.
The scaling multiplies the first derivatives by a factor of $1/\mathop{\rm side}(\omega)$, and the second derivatives by a factor of $1/\mathop{\rm side}(\omega)^2$.

If $\omega'$ is a scaled cube containing $x$, then
\begin{align*}
    d(x, Q^c) &\ge d(\omega', Q^c)
        \\&\ge d(\omega, Q^c) - {\sqrt{d} \over 2} {1 \over \sqrt{d}} \mathop{\rm side}(\omega)
        \\&\ge (1/2) \mathop{\rm side}(\omega),
\end{align*}
so the original cube $\omega$ has side length at most $2 d(x, Q^c)$.

In the other direction,
\begin{align*}
    d(x, Q^c) &\le d(\omega', Q^c) + \mathop{\rm diam}(\omega')
        \\&\le d(\omega, Q^c) + \frac32 \sqrt{d} \mathop{\rm side}(\omega')
        \\&\le (8 + 3\sqrt{d}/2) \mathop{\rm side}(\omega),
\end{align*}
so the original cube has side at least $\beta d(x,Q^c)$ where $\beta = 1/(8 + 3\sqrt{d}/2)$.

There are at most $2^d$ scaled cubes of a certain side length containing any point $x$. Let $N = \lfloor \log_2(2/\beta) + 1 \rfloor$. Then the total number of cubes containing $x$ of any side length is at most $N2^d$.

We now construct the smooth function as follows:
$$\Delta_Q(x) := {1 \over \beta} \sum_{\omega} \mathop{\rm side}(\omega) h_\omega(x).$$
This sum is locally finite, and every point $x \in Q$ is contained in at least one unscaled cube $\omega$ with $\mathop{\rm side}(\omega) \ge Bd(x, Q^c)$, so $\Delta_Q(x) \ge d(x, Q^c)$.
We also have the upper bounds that we want: $$\Delta_Q(x) \le {2 \over \beta} \sum_{\omega} d(x,Q^c) h_\omega(x) \le {N2^{d+1} \over \beta} d(x, Q^c).$$ Set $C_1 := N2^{d+1}/\beta$. For the first derivative, $|\nabla h_\omega| \le \max |\nabla h| /\mathop{\rm side}(\omega)$, so $$|\nabla \Delta_Q(x)| = \left|\frac1\beta \sum_\omega \mathop{\rm side}(\omega) \nabla h_\omega\right| \le C_2,$$ where $C_2 = N2^d \beta^{-1} \max |\nabla h(x)|$. In the same way we get the bound$$\left|{\partial^2 \Delta_Q \over \partial x_a \partial x_b}\right| \le {1 \over \beta} \sum_\omega {1 \over \mathop{\rm side}(\omega)} |\partial_a \partial_b h_\omega(x)| \le {C_3 \over d(x,Q^c)},$$ where the constant is $C_2 = N2^d\beta^{-2}\max |\partial_a \partial_b h(x)|$. These are the bounds that we want, with $C = \max\{C_1, C_2, C_3\}$. \qed

\begin{lemma}[Hedberg, \cite{hedbergappx}, Lemma 4]\label{lemmafour} Let $\delta = \delta(x)$ be $d(x, Q^c)$.
If $Q$ is a bounded open set, there is a sequence $h_j \in C^\infty_c(Q)$ with $0 \le h_j \le 1$ everywhere, $h_j(x)=1$ if $\delta(x) > 1/j$, and~\begin{align*}\left|\nabla h_j(x)\right| &\le {\xi(\delta(x))\over{j}},&\left|{\partial^2 \over \partial x_a \partial x_b} h_j(x)\right| &\le{2\xi(\delta(x)) \over j\delta}.\end{align*}\end{lemma}

\angry Proof.
Let $\Delta_Q(x)$ be the smooth distance function from Stein \S VI.2 Theorem 2 \cite{steinsingular}, quoted below as Lemma~\ref{partitionofunitygivesdistancefunction}. Let $C$ be the constant.

Set ${\rm H}_j(t) := \int_0^t \eta_{1/Cj}(\tau) \,d\tau$,
where $\eta_m$ with $m = 1/Cj$ is obtained from Lemma \ref{somerealvariablenonsense}.
Our sequence is $h_j(x) := {\rm H}_j(\Delta_Q(x))$.
These are compactly supported because $\Delta_Q(x) \le C\delta(x)$ and ${\rm H}_j(t)$ is zero for small enough $t$. We need to prove the bounds on the first and second derivatives.

The gradients of this sequence of functions are $\nabla h_j = \eta_j(\Delta_Q) \nabla \Delta_Q$, and
$${\partial^2 h_j \over \partial x_a \partial x_b} = \eta_j(\Delta_Q) {\partial^2 \Delta_Q \over \partial x_a \partial x_b} + \eta'_j(\Delta_Q) {\partial \Delta_Q \over \partial x_a} {\partial \Delta_Q \over \partial x_b}.$$
We plug in the bounds $|\eta_{1/Cj}(t)| \le \xi(t) / Cj$, $|\eta'_{1/Cj}(t)| \le |\xi(t)|/Cjt$ from the first lemma above, and the bounds $\delta(x) \le \Delta_Q(x)$ and $|\nabla\Delta_Q| \le C$ and $|\partial_a\partial_b\Delta_Q| \le C/\delta(x)$ from the second lemma above.

The resulting bounds are $|\nabla h_j| \le \xi(\delta) / j$ and $|\partial_a \partial_b h_j|\le \xi(\delta)/j \delta$, and this is what we want.
\qed

\end{subsection}
\end{section}

\bibliographystyle{acm}
\bibliography{note}
\end{document}